\documentclass[12pt]{amsart} 
\usepackage{amscd}
\usepackage{epsfig, amssymb}
\allowdisplaybreaks

\newcommand{\Cc}{\operatorname{CC}}
\newcommand{\ct}{\widetilde{C}}
\newcommand{\Cend}{\operatorname{Cend}}
\newcommand{\Chom}{\operatorname{Chom}}
\newcommand{\coh}{cohomology}
\newcommand{\Conf}{\operatorname{Conf}}
\newcommand{\Cur}{\operatorname{Cur}}
\newcommand{\Hom}{\operatorname{Hom}}
\renewcommand{\H}{\operatorname{H}}
\newcommand{\HC}{\operatorname{HC}}
\newcommand{\id}{\operatorname{id}}
\newcommand{\Ker}{\operatorname{Ker}}

\newcommand{\E}{\mathcal{E}}
\newcommand{\F}{\mathcal{F}}
\newcommand{\Res}{\operatorname{Res}}
\newcommand{\sign}{\operatorname{sign}}
\newcommand{\vir}{\operatorname{Vir}}

\newcommand{\gtg}{\mathfrak{g}}

\newcommand{\na}{\mathbb{A}}
\newcommand{\nc}{\mathbb{C}}
\newcommand{\nz}{{\mathbb{Z}}}

\newtheorem{thm}{Theorem}[section]
\newtheorem{lm}{Lemma}[section]
\newtheorem{prop}{Proposition}[section]
\newtheorem{crl}{Corollary}[section]

\theoremstyle{definition}

\newtheorem{df}{Definition}[section]
\newtheorem{ex}{Example}[section]

\theoremstyle{remark}
\newtheorem{remark}{Remark}[section]

\newtheorem*{ack}{Acknowledgment}

\newcommand{\CC}{\mathbb C}

\newcommand{\ZZ}{\mathbb Z}

\numberwithin{equation}{section}
\renewcommand\vec[1]{{\boldsymbol{#1}}}
\def\i{{\mathrm{i}}}       % imaginary unit
\def\Cset{\mathbb{C}}       % complex numbers
\def\Zset{\mathbb{Z}}       % integers
\def\al{\alpha}                         %%% some abbreviations
\def\be{\beta}
\def\ga{\gamma}
\def\Ga{\Gamma}
\def\de{\delta}
\def\De{\Delta}
\def\ep{\varepsilon}
\def\la{\lambda}
\def\La{\Lambda}

\def\th{\theta}

\def\g{{\mathfrak{g}}}      % Lie algebra g
\def\h{{\mathfrak{h}}}      % Lie algebra h
\def\n{{\mathfrak{n}}}      % Lie algebra n
\def\u{{\mathfrak{u}}}      % Lie algebra u
\def\sl{{\mathfrak{sl}}}
\def\what{\widehat}
\def\d{\partial}
\def\tt{\otimes}
\def\ov{\overline}
\def\Vec{{\mathcal{V}}ec}
\def\Vect{{\mathcal{V}}ect}
\def\M{{\mathcal{M}}}

\def\symm{S}
\def\surjto{\twoheadrightarrow}                %%  -->> surjection
\def\wti{\widetilde}
\DeclareMathOperator{\rank}{rank}

\DeclareMathOperator{\Sym}{S}       % symmetric power
\DeclareMathOperator{\Spec}{Spec}
\DeclareMathOperator{\Der}{Der}

\DeclareMathOperator{\Lie}{Lie}
\begin{document}

\title
{Cohomology of Conformal Algebras}

\author[B. Bakalov]{Bojko Bakalov}
\author[V. G. Kac]{Victor G. Kac}
\author[A. A. Voronov]{Alexander A. Voronov}
\address{Department of Mathematics, M.I.T., 77 Massachusetts Ave.,
Cambridge, MA 02139-4307}
%\curraddr{}
\email{bakalov@math.mit.edu}
\email{kac@math.mit.edu}
\email{voronov@math.mit.edu}
\thanks{Research of Bakalov and Kac was supported in part by 
NSF grant \#DMS-9622870 and
of Voronov by an AMS Centennial Fellowship.}

%\subjclass{Primary 16S80, 16E40; Secondary 81S10, 55P62}
\date{August 6, 1998}
\dedicatory{To Bertram Kostant on his seventieth birthday}

%\begin{abstract}
%\end{abstract}

\maketitle
\tableofcontents
%%%%%%%%%%%%%%%%%%%%%%%%%%%%%%%%%%%%%%%%%%%%%%%%%%%%%%

\section*{Introduction}
The notion of a conformal algebra encodes an axiomatic description of
the operator product expansion of chiral fields in conformal field
theory.  On the other hand, it is an adequate tool for the study of
infinite-dimensional Lie algebras satisfying the locality property
\cite{K}--\cite{K2}, \cite{DK}.  Likewise, conformal modules over a
conformal algebra $A$ correspond to conformal modules over the
associated Lie algebra $\Lie A$ \cite{CK}.  The main examples of Lie
algebras $\Lie A$ are the Lie algebras ``based'' on the punctured
complex plane $\nc^\times$, namely the Lie algebra $\Vect \nc^\times$ of vector
fields on $\nc^\times$ ($=$ Virasoro algebra) and the Lie algebra of maps
of $\nc^\times$ to a finite-dimensional Lie algebra ($=$ loop algebra).  Their
irreducible conformal modules are the spaces of densities on $\nc^\times$
and loop modules, respectively, \cite{CK}.  Since complete reducibility
does not hold in this case (cf.\ \cite{F},
\cite{CKW}), one may expect that their cohomology theory is very
interesting.

In the present paper we develop a cohomology theory of conformal
algebras with coefficients in an arbitrary module.  We introduce the
basic and the reduced complexes, the latter being a quotient of the
former.  The basic complex turns out to be isomorphic to the Lie
algebra complex for the so-called annihilation subalgebra $(\Lie A)_-$
of $\Lie A$.  For the main examples the annihilation subalgebra turns
out to be its complex-plane counterpart (\emph{i.e.}, $\nc^\times$ 
is replaced by
$\nc$).  The cohomology of these Lie algebras has been extensively
studied in \cite{GF,GF2,FF,Fe,F,Fe2}.  This allows us to compute the
cohomology of the conformal algebra $A$, which in its turn captures
main features of the cohomology of the Lie algebra $\Lie A$.  As a
byproduct of our considerations, we compute the \coh\ of a current Lie
algebra on $\nc$ with values in an irreducible highest-weight module
(see Theorem~\ref{cohcurnontr}), which has been known only when the
module is trivial \cite{Fe}.

The first \coh\ theory in the context of operator product expansion
was the \coh\ theory of vertex algebras and conformal field
theories introduced in \cite{KV}. The \coh\ theory of the present paper
relates to the \coh\ theory of \cite{KV} as much as
Chevalley--Eilenberg \coh\ of Lie algebras relates to Hochschild (or
more exactly, Harrison) \coh\ of commutative associative algebras. The
two theories possess standard properties of
\coh\ theories. For example, the \coh\ of \cite{KV} describes
deformations of vertex algebras, and the \coh\ of this paper
describes same of conformal algebras. However, the \coh\ of
\cite{KV} is hard to compute, whereas this paper offers 
the computation of \coh\ in most of the important examples.

The paper is organized as follows.  In Section~\ref{s1} we recall the
definition of a conformal algebra and of a (conformal) module 
over it and describe
their relation to formal distribution Lie algebras and conformal
modules.

In Section~\ref{s2} we construct the \emph{basic} complex
$\widetilde{C}^\bullet(A,M)$ and its quotient, the \emph{reduced}
complex $C^\bullet(A,M)$, for a module $M$ over a conformal algebra
$A$. These complexes define the basic and reduced
\coh\ of a conformal algebra $A$.

In Section~\ref{s4} we show that this cohomology parameterizes
$A$-module extensions, abelian conformal-algebra extensions,
first-order deformations, etc. (Theorem~\ref{exts}).

In Section~\ref{s5} we construct the dual, homology
complexes.  In Section~\ref{s7} we define the exterior
multiplication, contraction and module structure for the basic
complex.

In Section~\ref{s8} we prove that the basic complex is isomorphic to
the Lie algebra complex of the annihilation algebra (Theorem
6.1). Along with Proposition~\ref{prop:1.1} this implies, in
particular, that basic cohomology can be defined via a derived
functor. Apparently this is not the case for the reduced complex.

In Section~\ref{s9} we compute the cohomology with trivial
coefficients of the
Virasoro conformal algebra $\vir$ both for the basic and reduced complexes
(Theorem 7.1).  As one
could expect, the calculation and the result are closely related to
Gelfand--Fuchs's calculation of the cohomology of $\Vect \nc^\times$
\cite{GF}. We also compute both cohomologies  of $\vir$ with coefficients
in the modules of densities (Theorem 7.2). This result is closely related
to the work of Feigin and Fuchs \cite {FF, F}.

In Section~\ref{scc} we compute the cohomology of the current
conformal algebras both with trivial coefficients (Theorem 8.1) and
with coefficients in current modules (Theorem 8.2). This allows us,
in particular, to classify abelian extensions of current algebras
(Remark 8.1). Of course, abelian extensions of $\vir$ can be
classified by making use of Theorem 7.2.  This problem has been solved
earlier by M. Wakimoto and one of the authors of the present paper by
a lengthy but direct calculation; however, in the case of current
algebras the direct calculation is all but impossible.

In Section~\ref{shc} we briefly discuss the analogues of Hochschild
and cyclic cohomology for associative conformal algebras and of
Leibniz \coh.

In Section~\ref{s12} we indicate how to generalize our cohomology
theory to the case of conformal algebras in several indeterminates and
discuss its relation to cohomology of Cartan's filtered Lie algebras.

In Section~\ref{s11} we introduce anticommuting higher differentials
which may be useful for computing the cohomology of the basic
complex with non-trivial coefficients.

In Section~\ref{s13} we briefly discuss the relation
of our cohomology theory to $\Lie$ algebras in a general
pseudo-tensor category introduced in \cite{BD}.

In the last Section~\ref{sop} we list several open questions.

Unless otherwise specified, all vector spaces, linear maps and
tensor products are considered over the field $\CC$ of complex
numbers.  We will use the divided-powers notation $\la^{(m)} = \la^m / m!$, 
$m\in\ZZ_+$, where $\ZZ_+$ is the set of non-negative integers.

\begin{ack}
The second author is grateful to Jean-Louis Loday for inspiring
discussions on Leibniz algebras, whose conformal version, see
Section~\ref{leibniz}, seems to be an essential notion in the case of
nonlocal fields.
\end{ack}
%%%%%%%%%%%%%%%%%%%%%%%%%%%%%%%%%%%%%%%%%%%%%%%%%%%%%
\section{Preliminaries on conformal algebras and modules}\label{s1}
%%%%%%%%%%%%%%%%%%%%%%%%%%%%%%%%%%%%%%%%%%%%%%%%%%%%%

\begin{df}
  \label{df:1.1}
A (\emph{Lie$)$ conformal algebra} is a $\CC [\partial ]$-module $A$
endowed with a $\lambda$-bracket $[a_{\lambda} b]$ which defines a
linear map $A \otimes A \to A [ \lambda ]$, where $A [ \lambda ] = \nc
[ \lambda ] \tt A$, subject to the following axioms:

\begin{description}

\begin{sloppypar}
\item[Conformal sesquilinearity] $[ \partial a_{\lambda} b] =
  - \lambda [a_{\lambda} b]$, $[ a_{\lambda} \partial b ] \linebreak[1] =
  \linebreak[0] (\partial + \lambda ) [a_{\lambda} b] $;
\end{sloppypar}

\item[Skew-symmetry]  $ [a_{\lambda} b] = -[b_{- \lambda - \partial
    } a]$;

\item[Jacobi identity]  $[ a_{\lambda} [b_{\mu} c]] =
        [[a_{\lambda} b]_{\lambda + \mu} c] +
        [b_{\mu} [ a_{\lambda} c]]$.

\end{description}

\end{df}

Conformal algebras appear naturally in the context of formal
distribution Lie algebras as follows.  Let $\gtg$ be a vector
space.  A $\gtg$-\emph{valued formal distribution} is a series of
the form $a(z) = \sum_{n \in \ZZ} a_n z^{-n-1}$, where $a_n \in
\gtg$ and $z$ is an indeterminate.  We denote the space of such
distributions by $\gtg [[z, z^{-1}]]$ and the operator $\partial_z$
on this space by $\partial$.

Let $\gtg$ be a Lie algebra.  Two $\gtg$-valued formal
distributions are called \emph{local} if
\begin{displaymath}
  (z-w)^N [a(z), b(w)]= 0 \hbox{ for } N \gg 0 \, .
\end{displaymath}
This is equivalent to saying that one has an expansion of the form
\cite{K}:

\begin{align}
  \label{eq:1.1}
  [a(z) , b(w)] &= \sum^{N-1}_{j=0} \bigl(a(w)_{(j)} b(w) \bigr)
    \partial^{(j)}_w \delta (z-w),
\\
\intertext{where}
  \label{eq:1.2}
  a(w)_{(j)} b(w) &= \Res_z (z-w)^j [a(z) , b(w)] 
\\
\intertext{and}
  \notag  \delta (z-w) &= \sum_{n \in \ZZ}
         z^{-n-1} w^n .
\end{align}
Let $\F$ be a family of pairwise local $\gtg$-valued formal
distributions such that the coefficients of all distributions
from $\F$ span $\gtg$.  Then the pair $(\gtg , \F$) is called
a \emph{formal distribution Lie algebra}.

Let $\overline{\F}$ denote the minimal subspace of $\gtg [[z,z^{-1}]]$
containing $\F$ which is closed under all $j$-th products
(\ref{eq:1.2}) and $\partial$-invariant.  One knows that
$\overline{\F}$ still consists of pairwise local distributions
\cite{K}.  Letting
\begin{displaymath}
  [a_{\lambda} b] = \sum_{n \in \ZZ_+} \lambda^{(n)}
  a_{(n)} b ,
\end{displaymath}
%where $\lambda^{(n)} = \la^n/n!$, 
one endows $\overline{\F}$ with the structure of a conformal
algebra, which is denoted by $\Conf (\gtg , \F)$ \cite{DK, K}.

Conversely, given a conformal algebra $A$, one associates to it
the \emph{maximal} formal distribution Lie algebra ($\Lie A,A)$
as follows.

Let $\Lie A = A[t,t^{-1}] / (\partial + \partial_t) A [t,t^{-1}]$ and
let $a_n$ denote the image of $at^n$ in $\Lie A$.  Then the formula
$(a,b \in A, m, n \in
\ZZ)$:
\begin{equation}
  \label{eq:1.3}
  [a_m, b_n ] = \sum_{j \in \ZZ_+} \binom{m}{j} (a_{(j)}b)_{m+n-j}
\end{equation}
gives a well defined bracket making $\Lie A$ a Lie algebra. It forms a
formal distribution Lie algebra with the family of pairwise local
distributions $\F = \left\{ a(z) = \sum_{n \in \ZZ} a_n z^{-n-1}
\right\}_{a \in A}$.  We have: $\Conf (\Lie A , \F)
\simeq A$ via the map $a \mapsto a(z)$ \cite{K}.

The Lie algebra $\Lie A$ carries a derivation $T$ induced by $-
\partial_t$:
\begin{equation}
  \label{eq:1.4}
  T (a_n) = -na_{n-1} \, .
\end{equation}
It is clear from (\ref{eq:1.3}) that the $\CC$-span of the $a_n$ with
$n \in \ZZ_+, a \in A$, is a $T$-invariant subalgebra of the Lie
algebra $\Lie A$.  This subalgebra is denoted by $(\Lie A)_-$ and is
called the \emph{annihilation Lie algebra} of $A$.  The semidirect sum
$(\Lie A)^- = \CC T + (\Lie A)_-$ is called the \emph{extended
annihilation Lie algebra}.

If one drops the skew-symmetry in the definition of a Lie algebra
$\gtg$, but keeps the Leibniz version of the Jacobi identity
$[a,[b,c]] = [[a,b] ,c] + [b, [a,c]]$, then $\gtg$ is called a (left)
\emph{Leibniz algebra}, see \cite{loday}. If one also drops the
condition of locality on $\F$, then $(\gtg, \F)$ is called a
\emph{formal distribution Leibniz algebra}. In this case $\Conf (\gtg ,
\F)$ is a \emph{Leibniz conformal algebra}, \emph{i.e}., the
skew-symmetry axiom in the definition of a Lie conformal algebra is
dropped.

\begin{df}
  \label{df:1.2}
A \emph{module} $M$ over a Lie conformal algebra $A$ is a $\CC
[\partial]$-module endowed with the $\lambda$-action
$a_{\lambda}v$ which defines a map $A \otimes M \to M [[
\lambda ]]$ such that
\begin{align}
\label{module}
 & a_{\lambda} (b_{\mu} v) - b_{\mu} (a_{\lambda}v) =
[a_{\lambda}b]_{\lambda + \mu} v ,
\\
& (\d a)_\la v = -\la a_\la v, \quad
a_\la (\d v) = (\d+\la) a_\la v .
\end{align}
If $a_{\lambda}v \in M [\lambda]$ for all $a \in A$, $v
\in M$, then the $A$-module $M$ is called \emph{conformal}. If $M$ is
finitely generated over $\nc[\d]$, $M$ is simply called \emph{finite}.
\end{df}

\begin{df}[\cite{DK}]
\label{clm}
A \emph{conformal linear map} from an $A$-module $M$ to an $A$-module
$N$ is a $\nc$-linear map $f\colon  M \to N[\lambda]$, 
denoted $f_\lambda\colon  M\to N$, 
such that $f_\lambda \d = (\d+\la) f_\lambda$. The space
of such maps is denoted $\Chom (M,N)$. It has canonical structures
of a $\nc[\d]$- and an $A$-module:
\begin{align*}
&(\d f)_\la  = - \la f_\la,
\\
&(a_\mu f)_\la m = a_\mu (f_{\la - \mu} m) - f_{\la-\mu}(a_\mu m),
\end{align*}
where $a \in A$, $m \in M$, and $f \in \Chom(M,N)$. When the two
modules $M$ and $N$ are conformal and finite, the module $\Chom(M,N)$
will also be conformal.
\end{df}

For a finite module $M$, let $\Cend M = \Chom(M,M)$ denote the space
of conformal linear endomorphisms of $M$. Besides the $A$-module
structure, $\Cend M$ carries the natural structure
\[
(f_\la g)_\mu m = f_\la (g_{\mu - \la} m), \qquad f, g \in \Cend M, m
\in M,
\]
of an associative conformal algebra in the sense of the following
definition, see \cite{K2}.

\begin{df}
\label{df:1.4}
An \emph{associative conformal algebra} is a $\CC [\partial ]$-module
$A$ endowed with a $\lambda$-multiplication $a_{\lambda} b$ which
defines a linear map $A \nolinebreak \otimes \nolinebreak A \to
A[\lambda]$ subject to the following axioms:
\begin{description}

\item[Conformal sesquilinearity] $(\partial a)_{\lambda} b =
  - \lambda a_{\lambda} b , a_{\lambda} \partial b =
  (\partial + \lambda ) a_{\lambda} b $;

\item[Associativity]  $a_{\lambda} (b_{\mu} c) =
        (a_{\lambda} b)_{\lambda + \mu} c$.

\end{description}
\end{df}

The $\la$-bracket $[a_\la b ] = a_\la b - b_{-\la - \d} a$ makes an
associative conformal algebra, in particular, $\Cend M$, a Lie
conformal algebra. $\Cend M$ with this structure is denoted
$\operatorname{gc} M$ and called the \emph{general Lie conformal
algebra of a module} $M$ \cite{DK,K2}.

Given an associative conformal algebra $A$, a \emph{left\/} 
(or \emph{right}) \emph{module} 
$M$ over it may be defined naturally, for example,
like in Definition~\ref{df:1.2}. A \emph{bimodule} may be defined by
adding the axiom $a_{\lambda} (m_{\mu} b) = (a_{\lambda}m)_{\lambda +
\mu} b$ to the list of those for a left and right module. A (bi)module
is called \emph{conformal}, provided the action(s) satisfy the usual
polynomiality conditions. The structure of a conformal bimodule on $M$
is equivalent to an extension of the associative conformal algebra
structure to the space $A \oplus \epsilon M$, where $\epsilon^2 =0$.

We will be working with Lie conformal algebras and modules over them
throughout the paper, except when we discuss Hochschild \coh\ in
Section~\ref{hochschild}. We will therefore usually shorten the term
``Lie conformal algebra'' to ``conformal algebra''.

Conformal modules over conformal algebras appear naturally in the
context of conformal modules over formal distribution Lie algebras as
follows.  Let $(\gtg , \F)$ be a formal distribution Lie algebra and
let $V$ be a $\gtg$-module.  Suppose that $\E$ is a family of
$V$-valued formal distributions which spans $V$ and such that any
$a(z) \in \F$ and $v(z) \in \E$ form a local pair, \emph{i.e.},
\begin{displaymath}
  (z-w)^N a(z) v(w) =0 \quad \hbox{for } N \gg 0 \, .
\end{displaymath}
Then $(V,\E)$ is called a \emph{conformal} $(\gtg , \F)$-module.
As before, we have:
\begin{align}
  \label{eq:1.5}
  a(z) v(w) &= \sum^{N-1}_{j=0}
     \bigl( a(w)_{(j)} v(w) \bigr) \partial^{(j)}_w \delta
     (z-w), 
\\
\intertext{where}
  \label{eq:1.6}
  a(w)_{(j)}v (w) &= \Res_z (z-w)^j a(z) v(w) \, .
\end{align}
Let $\overline{\E}$ denote the minimal subspace of $V[[z,z^{-1}]]$
containing $\E$ which is closed under all $j$-th actions
(\ref{eq:1.6}) and is $\partial$-invariant. One knows that all
pairs $a(z) \in \overline{\F}$ and $v(z) \in \overline{\E}$ are
still local \cite{K, K2}.  Letting
\begin{displaymath}
  a_{\lambda} v = \sum_{n \in \ZZ_+} 
  \lambda^{(n)} a_{(n)} v ,
\end{displaymath}
%where $\lambda^{(n)} = \la^n/n!$, 
one endows $\overline{\E}$ with the structure of a conformal
$\overline{\F}$-module \cite{K, K2}.

Conversely, given a conformal $A$-module $M$, one associates to it the
\emph{maximal} conformal $(\Lie A,A)$-module $(V(M),M)$ in a way
similar to the one the Lie algebra $\Lie A$ has been constructed.  We
let $V(M) =M [t,t^{-1}] / (\partial + \partial_t) M [t,t^{-1}]$, with
the well-defined $\Lie A$-action
\begin{equation}
  \label{eq:1.7}
  a_m v_n = \sum_{j \in \ZZ_+} \binom{m}{j}
  (a_{(j)} v)_{m+n-j} \, ,
\end{equation}
where, as before, $v_n$ stands for the image of $vt^n$ in $V(M)$
\cite{K}.  As before, we denote by $V(M)_-$ the $\CC$-span of the
$v_n$, where $v \in M$, $n \in \ZZ_+$.  It is clear from
(\ref{eq:1.7}) that $V(M)_-$ is a $(\Lie A)^-$- and a $(\Lie
A)_-$-submodule of $V(M)$.

The following obvious observation plays a key role in
representation theory of conformal algebras \cite{CK}.

\begin{prop}
  \label{prop:1.1}
A module $M$ over a conformal algebra $A$ 
% is the same as
carries the natural structure of a module over the extended
annihilation Lie algebra $(\Lie A)^-$. This correspondence establishes
an equivalence of the category of $A$-modules and that of $(\Lie
A)^-$-modules. The $A$-module $M$ is conformal, iff as a $(\Lie
A)^-$-module it satisfies the condition
\begin{equation}
  \label{eq:1.8}
  a_nv =0 \quad \hbox{for } a \in A , \quad v \in V , 
  \quad n \gg 0 \, .
\end{equation}

\end{prop}

\begin{remark}
  \label{rem:1.1}
As a $(\Lie A)^-$-module, a
conformal $A$-module $M$ is isomorphic to the module
$V(M) / V(M)_-$.

\end{remark}

\begin{remark}
  \label{rem:1.2}
\cite{K}. One can show that the map $A \mapsto (\Lie A,A)$
(respectively, $M \mapsto (V(M),M)$) establishes a bijection between
isomorphism classes of conformal algebras (respectively, of conformal
modules over conformal algebras) and equivalence classes of formal
distribution Lie algebras $(\gtg , \F)$ (respectively, of conformal
modules over $(\Lie A,A)$).  By definition, all formal distribution
Lie algebras $((\Lie A) / I, \F)$, where $I$ is an ideal of $\Lie A$
having trivial intersection with $A$, and $\overline{\F} =A$ are
\emph{equivalent} (and similarly for modules).

\end{remark}

\begin{ex}
  \label{ex:1.1}
Let $\gtg$ be a Lie algebra and let $\tilde{\gtg} = \gtg
[t,t^{-1}]$ be the associated loop ($=$ current) algebra (with
the obvious bracket:  $[at^m , bt^n ]= [a,b] t^{m+n}$, $a,b \in
\gtg$, $m,n \in \ZZ$).  For $a \in \gtg$ let $a(z) = \sum_{m \in
  \ZZ}(at^m) z^{-m-1} \in \tilde{\gtg} [[ z, z^{-1}]]$.  Then
\begin{displaymath}
  [a(z), b(w)] = [a,b] (w) \delta (z-w) \, ,
\end{displaymath}
hence the family $\F = \left\{ a(z) | a \in \gtg \right\}$
consists of pairwise local formal distributions and
$(\tilde{\gtg}, \F)$ is a formal distribution Lie algebra.  Note that
\begin{displaymath}
  \overline{\F} = \CC [ \partial ] \F \simeq \CC [ \partial ]
  \otimes \gtg
\end{displaymath}
is a conformal algebra with the $\lambda$-bracket
\begin{displaymath}
  [a_{\lambda} b] = [a,b] \, , \quad a,b \in \gtg \, .
\end{displaymath}

\end{ex}
This conformal algebra is called the \emph{current conformal
algebra} associated to $\gtg$ and is denoted by $\Cur \gtg$.
Note that $\Lie \, (\Cur \gtg , \F) \simeq \tilde{\gtg}$, hence
$\tilde{\gtg}$ is the maximal formal distribution algebra.  The
corresponding annihilation algebra is $\tilde{\gtg}_- = \gtg [t]$
and the extended annihilation algebra is $\CC \partial_t + \gtg
[t]$.

Given a $\gtg$-module $U$, one may associate the conformal
$\tilde{\gtg}$-module $\tilde{U} = U [t, t^{-1}]$ with the obvious action
of $\tilde{\gtg}$, and the conformal $\Cur \gtg$-module $M_U =
\CC [\partial ] \otimes U$ defined by
\begin{displaymath}
  a_{\lambda} u = au \, , \quad a \in \gtg \, , \, u \in U \, .
\end{displaymath}
We have:  $V(M_U) \simeq \tilde{U}$ as $\tilde{\gtg}$-modules.

It is known that, provided that $\gtg$ is finite-dimensional
semisimple, the $\Cur \gtg$-modules $M_U$, where $U$ is a
finite-dimensional irreducible $\gtg$-module, exhaust all finite
irreducible non-trivial $\Cur \gtg$-modules \cite{CK}.

\begin{ex}\label{ex:1.2}
Let $\Vect \CC^\times$ denote the Lie algebra of all regular vector
fields on $\CC^\times$.  The vector fields $t^n \partial_t$ $(n \in \ZZ)$
form a basis of $\Vect \CC^\times$ and the formal distribution $L(z)=
- \sum_{n \in \ZZ} (t^n \partial_t) z^{-n-1}$ is local (with respect
to itself), since
\begin{displaymath}
  [L(z), L(w)] = \partial_w L(w) \delta (z-w) +
     2L (w) \delta'_w (z-w) \, .
\end{displaymath}
 Hence $(\Vect \, \CC^\times , \left\{ L \right\})$ is a formal
 distribution Lie algebra.  The associated conformal algebra
\begin{displaymath}
\vir = \CC [ \partial ] L \, , \quad [L_{\lambda}L] =
( \partial + 2 \lambda )L
\end{displaymath}
is called the \emph{Virasoro conformal algebra}.

Note that $\Lie \, (\vir \, , \left\{ L \right\}) \simeq \Vect \,
\CC^\times$, hence $\Vect \, \CC^\times$ is the maximal formal distribution
algebra.  The corresponding annihilation algebra $(\Vect \,
\CC^\times)_- = \Vect \, \CC$, the Lie algebra of regular vector
fields on $\CC$, and $(\Vect \, \CC^\times)^-$ is isomorphic to the
direct sum of $(\Vect \, \CC^\times)_-$ and the $1$-dimensional Lie algebra.

It is known that all free non-trivial $\vir$-modules of rank $1$ over
$\CC [ \partial ]$ are the following ones $( \Delta ,
\alpha \in \CC)$:
\begin{displaymath}
  M_{\Delta , \alpha} = \CC [ \partial ]v \, , \quad
  L_\la v= ( \partial + \alpha + \Delta \lambda )v \, .
\end{displaymath}
We have: $V(M_{\Delta , \alpha}) \simeq \CC [ t,t^{-1}] e^{- \alpha t}
(dt)^{1- \Delta}$ as $\Vect \, \CC$-modules.  The module $M_{\Delta ,
\alpha}$ is irreducible, iff $ \Delta \neq 0$.  The module $M_{0,
\alpha}$ contains a unique non-trivial submodule $(\partial + \alpha)
M_{0, \alpha}$ isomorphic to $M_{1, \alpha}$.  It is known that the
modules $M_{\Delta , \alpha}$ with $\Delta \neq 0$ exhaust all finite
irreducible non-trivial $\vir$-modules \cite{CK}.

\end{ex}

It is known \cite{DK} that the conformal algebras $\Cur \gtg$, where
$\gtg$ is a finite-dimensional simple Lie algebra, and $\vir$ exhaust
all finite simple conformal algebras.  For that reason we shall
discuss mainly these two examples in what follows.

%%%%%%%%%%%%%%%%%%%%%%%%%%%%%%%%%%%%%%%%%%%%%%%%%%%%%%%%%%%%%%%
\section{Basic definitions}\label{s2}
%%%%%%%%%%%%%%%%%%%%%%%%%%%%%%%%%%%%%%%%%%%%%%%%%%%%%%%%%%%%%%%

%\subsection{Cochains}

\begin{df}
\label{cochains}
An \emph{$n$-cochain $(n \in \nz_+)$ of a conformal algebra $A$ with coefficients in a
module $M$ over it} is a $\nc$-linear map
\begin{align*}
\gamma: A^{\otimes n} & \to M[\lambda_1, \dots, \lambda_n]\\
a_1 \otimes \dots \otimes a_n & \mapsto \gamma_{\lambda_1, \dots,
\lambda_n}(a_1, \dots, a_n),
\end{align*}
where $M[\lambda_1, \dots, \lambda_n]$ denotes the space of
polynomials with coefficients in $M$, satisfying the following
conditions:
\begin{description}
%\begin{sloppypar}
%\item for any $a_1, \dots, a_n$ in $A$,
%$\alpha_{\lambda_1, \dots, \lambda_n} (a_1, \dots, a_n)$ is polynomial
%in $\lambda_1, \linebreak[0] \dots, \linebreak[1] \lambda_n$;
%\end{sloppypar}
\item[Conformal antilinearity] $\gamma_{\lambda_1, \dots, \lambda_n}
(a_1, \dots, \partial a_i, \dots, a_n)$ \\  $= - \lambda_i
\gamma_{\lambda_1, \dots, \lambda_n} (a_1, \dots, a_i, \dots, a_n)$
for all $i$;
\item[Skew-symmetry] $\gamma$ is skew-symmetric with respect to simultaneous
permutations of $a_i$'s and $\lambda_i$'s.
%\item[Locality] $\gamma_{\lambda_1, \dots, \lambda_n}(a_1, \dots, a_n)$
%depends polynomially on the $\lambda_i$'s for any fixed $a_1, \dots,
%a_n \in A$.
\end{description}
We let $A^{\otimes 0} = \nc$, as usually, so that a 0-cochain $\gamma$
is an element of $M$.
Sometimes, when the module $M$ is not conformal, one may consider
formal power series instead of polynomials in this definition.
\end{df}

%\subsection{Differential}

We define a differential $d$ of a cochain $\gamma$ as follows:
\begin{multline*}
(d\gamma)_{\lambda_1, \dots, \lambda_{n+1}} (a_1,\dots,a_{n+1})
\\
\begin{split}
& = \sum_{i=1}^{n+1}  (-1)^{i+1} {a_i}_{\lambda_i}
\gamma_{\lambda_1, \dots, \widehat\lambda_i,
\dots, \lambda_{n+1}} (a_1,\dots, \widehat a_i,\dots, a_{n+1})
\\
& 
\begin{split}
+ \sum^{n+1}_{\substack{i,j = 1\\i < j}}
(-1)^{i+j} \gamma_{\lambda_i+\lambda_j,
\lambda_1, \dots, \widehat\lambda_i, \dots, \widehat\lambda_j, \dots,
\lambda_{n+1}} ([{a_i} _{\lambda_i} a_j], a_1, \dots, \widehat a_i,\dots, 
\widehat a_j,\\
\dots, a_{n+1}),
\end{split}
\end{split}
\end{multline*}
where $\gamma$ is extended linearly over the polynomials in
$\lambda_i$. In particular, if $\gamma \in M$ is a 0-cochain, then
$(d\gamma)_\lambda(a) = a_\lambda \gamma$.

\begin{remark}
\label{l+m}
Conformal antilinearity implies the following relation for an
$n$-cochain $\gamma$:
\[
\gamma_{\lambda+\mu, \lambda_1, \dots} ([a_\lambda b], a_1 , \dots )
= \gamma_{\lambda+\mu, \lambda_1, \dots} ([a_{-\partial-\mu} b ], a_1,
\dots ) .
\]
\end{remark}

\begin{lm}
\begin{enumerate}
\item The operator $d$ preserves the space of cochains;
\item $d^2 = 0$.
\end{enumerate}
\end{lm}

\begin{proof} 
1. The only non-trivial point in checking the skew-symmetry of $d\gamma$
amounts to the equation
\begin{equation*}
%\label{4}
\gamma_{\lambda + \mu, \lambda_1, \dots, \lambda_{n-1}} ([a _\lambda b],
a_1, \dots, a_{n-1})
=- \gamma_{\lambda + \mu, \lambda_1, \dots,
\lambda_{n-1}} ([b _\mu a], a_1, \dots, a_{n-1}),
\end{equation*}
which follows from Remark~\ref{l+m} and the skew-symmetry of
$[a_\lambda b]$.

2. To check that $d^2 = 0$, we will compute $d^2 \gamma$ for an
$n$-cochain $\gamma$.
\begin{multline*}
(d^2 \gamma)_{\lambda_1, \dots, \lambda_{n+2}} (a_1,\dots,a_{n+2})
\\
= \sum_{i=1}^{n+2}  (-1)^{i+1} {a_i} _{\lambda_i}
(d\gamma)_{\lambda_1, \dots, \widehat\lambda_i,
\dots, \lambda_{n+2}} (a_1,\dots,\widehat a_i,\dots, a_{n+2})
\\
+ \sum^{n+2}_{\substack{i,j = 1\\i < j}}
(-1)^{i+j} (d\gamma)_{\lambda_i+\lambda_j,
\lambda_1, \dots, \widehat\lambda_{i,j}, \dots,
\lambda_{n+2}} ([{a_i} _{\lambda_i} a_j], a_1, \dots, \widehat a_{i,j},\dots, 
a_{n+2})
\\
= \sum_{\substack{i,j=1\\i \ne j}}^{n+2} (-1)^{i+j}\sign\{j,i\}
{a_i}_{\lambda_i} \bigl( {a_j}_{\lambda_j} \gamma_{\lambda_1, \dots,
\widehat\lambda_{i,j},
\dots, \lambda_{n+2}} (a_1,
\dots,\widehat a_{i,j},\dots, a_{n+2}) \bigr)
\\
+ \sum_{\substack{i,j,k=1\\j<k, i \ne j,k}}^{n+2}  
(-1)^{i+j+k+1}\sign\{j,k,i\} {a_i} _{\lambda_i}
\gamma_{\lambda_j+\lambda_k, \lambda_1,  
\dots, \widehat\lambda_{i,j,k}, \dots, \lambda_{n+2}} ([{a_j}_{\lambda_j} a_k],
\\
a_1,\dots,\widehat a_{i,j,k}, \dots, a_{n+2})
\\
+ \sum^{n+2}_{\substack{i,j,k = 1\\i < j, k \ne i,j}}
(-1)^{i+j+k}\sign\{k,i,j\} {a_k}_{\lambda_k}
\gamma_{\lambda_i+\lambda_j, \lambda_1, \dots,
\widehat\lambda_{i,j,k}, \dots, \lambda_{n+2}} ([{a_i} _{\lambda_i}
a_j],
\\
a_1, \dots, \widehat a_{i,j,k}, \dots, a_{n+2})
\\
+ \sum^{n+2}_{\substack{i,j = 1\\i < j}} (-1)^{i+j}
[{a_i} _{\lambda_i}
a_j]_{\lambda_i+\lambda_j} \gamma_{\lambda_1, \dots,
\widehat\lambda_{i,j}, \dots, \lambda_{n+2}} (a_1,
\dots, \widehat a_{i,j}, \dots, a_{n+2})
\\
+ \sum^{n+2}_{\substack{\text{distinct }i,j,k,l = 1\\i < j, k<l}}
(-1)^{i+j+k+l}\sign\{i,j,k,l\}
\\
\times
\gamma_{\lambda_k+\lambda_l,
\lambda_i+\lambda_j,
\lambda_1, \dots, \widehat\lambda_{i,j,k,l}, \dots,
\lambda_{n+2}} ([{a_k} _{\lambda_k} a_l],
[{a_i} _{\lambda_i} a_j], a_1, \dots, \widehat a_{i,j,k,l},\dots,
a_{n+2})
\\
+ \sum^{n+2}_{\substack{i,j,k = 1\\i < j, k \ne i,j}}
(-1)^{i+j+k+1}\sign\{i,j,k\}
\\
\times
\gamma_{\lambda_i+\lambda_j+\lambda_k, \lambda_1, \dots,
\widehat\lambda_{i,j,k}, \dots, \lambda_{n+2}} ([[{a_i} _{\lambda_i}
a_j]_{\lambda_i+\lambda_j} a_k],
a_1, \dots, \widehat a_{i,j,k}, \dots, a_{n+2}),
\end{multline*}
where $\sign\{i_1,\dots,i_p\}$ is the sign of the permutation putting
the indices in the increasing order and $\widehat a_{i,j, \dots}$
means that $a_i, a_j, \dots$ are omitted.  Notice that each term in
the summation over $i,j,k,l$ is skew with respect to the permutation
$\begin{pmatrix} i & j & k & l\\k&l&i&j\end{pmatrix}$. Therefore, the
terms of that summation will cancel pairwise. The first and the forth
summations cancel each other, because $M$ is a conformal algebra
module:
\[
-{a_i}_{\lambda_i}({a_j}_{\lambda_j} m) +
{a_j}_{\lambda_j}({a_i}_{\lambda_i} m) +
[{a_i}_{\lambda_i}{a_j}]_{\lambda_i+\lambda_j} m = 0.
\] 
The second summation becomes equal to the third one after the
substitution $(ikj)$, except that they differ by a sign. Thus, they cancel
each other, as well. Finally, the sixth summation can be rewritten as
a summation over $i<j<k$ of the sum of three permutations of the initial
summand.  Precisely, in the first entry of $\gamma$, we will have
\begin{equation}
\label{jacobi1}
[[{a_i} _{\lambda_i}
a_j]_{\lambda_i+\lambda_j} a_k] - [[{a_i} _{\lambda_i}
a_k]_{\lambda_i+\lambda_k} a_j] + [[{a_j} _{\lambda_j}
a_k]_{\lambda_j+\lambda_k} a_i].
\end{equation}
Using Remark~\ref{l+m}, we can transform the sum \eqref{jacobi1}
inside $\gamma$ into
\[
[[{a_i} _{\lambda_i} a_j]_{\lambda_i+\lambda_j} a_k] - [[{a_i}
_{\lambda_i} a_k]_{-\partial - \lambda_j} a_j] + [[{a_j} _{\lambda_j}
a_k]_{-\partial-\lambda_i} a_i] ,
\]
which vanishes by the Jacobi identity and skew-symmetry in $A$. Thus,
we see that all of the terms in $d^2 \gamma$ cancel.
\end{proof}

%\subsection{Cohomology}

Thus the cochains of a conformal algebra $A$ with coefficients in a
module $M$ form a complex, which will be denoted
\[
\ct^\bullet =
\ct^\bullet (A,M)= \bigoplus_{n \in \nz_+} \ct^n (A,M).
\]
This complex is called the \emph{basic complex} for the $A$-module
$M$. This is not
yet the complex defining the right cohomology of a conformal algebra:
we need to consider a certain quotient complex.

Define the structure of a (left) $\nc[\partial]$-module on
$\ct^\bullet (A,M)$ by letting 
\begin{equation}\label{dwtic}
(\partial\cdot \gamma)_{\lambda_1, \dots, \lambda_n}(a_1, \dots, a_n) 
= \Bigl( \partial_M + \sum_{i=1}^n \lambda_i \Bigr) 
\gamma_{\lambda_1, \dots, \lambda_n}(a_1,\dots, a_n),
\end{equation}
where $\d_M$ denotes the action of $\d$ on $M$.

%In other words, we use the same
%$\Hom(A,\nc[[\lambda_1, \dots,
%\lambda_n]] \otimes_{\nc[\partial]} M)$ except that the tensor product
%is now taken over $\nc[\partial]$. 

\begin{lm}
\label{partial}
$d \partial = \partial d$, and therefore the graded subspace $\partial
\ct^\bullet \subset \ct^\bullet$ forms a subcomplex.
\end{lm}

\begin{proof}
The first summation in the differential transforms the factor
$\partial_M + \sum_{i=1}^n \lambda_i$ into $\partial_M +
\sum_{i=1}^{n+1} \lambda_i$, because of the conformal sesquilinearity of the
$\lambda$-bracket. The second summation does the same for more obvious
reasons.
\end{proof}

Define the quotient complex
\[
C^\bullet (A,M) = \ct^\bullet (A,M)/\partial \ct^\bullet (A,M) =
\bigoplus_{n \in \nz_+} C^n (A,M),
\]
called the \emph{reduced complex}.

\begin{df}

The \emph{basic \coh\ $\widetilde \H^\bullet(A,M)$ of a conformal
algebra $A$ with coefficients in a module $M$} is the \coh\ of the
basic complex $\ct^\bullet$.  The (\emph{reduced$)$ \coh}
$\H^\bullet(A,M)$ is the \coh\ of the \emph{reduced complex}
$C^\bullet = C^\bullet (A,M) = \ct^\bullet/\partial \ct^\bullet$.
%The (\emph{normalized$)$ \coh\ $H_0^\bullet(A,M)$} is the \coh\ of the
%subcomplex $C^\bullet_0 = \{c \in \ct^\bullet \; | \; \d c = 0 \}$.

\end{df}

\begin{remark}
The basic \coh\ $\widetilde \H^\bullet (A,M)$ is naturally a
$\nc[\d]$-mod\-ule, whereas the reduced \coh\ $\H^\bullet(A,M)$ is a complex
vector space.
\end{remark}

\begin{remark}
The exact sequence $0 \to \partial \ct^\bullet \to \ct^\bullet \to
C^\bullet \to 0$ gives the long exact sequence of \coh :
\begin{align}
\label{long-exact}
0 &\to \H^0(\partial \ct^\bullet) \to \widetilde \H^0(A,M) \to
\H^0(A,M) \to 
\\
\notag
&\to \H^1(\partial \ct^\bullet)
\to \widetilde \H^1 (A,M) \to \H^1(A,M) \to
\\
\notag
&\to  \H^2(\partial \ct^\bullet) \to
\widetilde \H^2 (A,M) \to \H^2(A,M) \to \dotsm
\end{align}
\end{remark}

\begin{prop}
\label{prop-dd}
In degrees $\ge 1$, the complexes $\ct^\bullet$ and $\d \ct^\bullet$ are
isomorphic under the map
\begin{align}
\label{dd}
\ct^\bullet \to \d \ct^\bullet, 
\qquad
\gamma  \mapsto \d \cdot \gamma.
\end{align}
Therefore, $\H^q(\d \ct^\bullet) \simeq \wti \H^q(A,M)$ for $q \ge 1$,
and the natural sequence $0 \to \Ker \d [0] \to \wti \H^0 (A,M) \to
\H^0 (\d \ct^\bullet) \to 0$, where $\Ker \d [0]$ is the subcomplex
$\Ker \d$ of $\ct^\bullet$, in fact concentrated in degree zero, is
exact.  When the module $M$ is $\nc[\d]$-free, the above isomorphisms
take place in all degrees $\ge 0$.
\end{prop}

\begin{proof}
Indeed, the modules $\ct^n(A,M)$, $n \ge 1$, are free over $\nc[\d]$,
because they are free over $\nc[\lambda_1]$. Lemma~\ref{partial} shows
that the map \eqref{dd} is a morphism of complexes. When $M$ is
$\nc[\d]$-free, this argument extends over to $n=0$.
\end{proof}

\begin{remark}
This proposition does not imply that in the long exact sequence
\eqref{long-exact}, the maps $\H^q (\d \ct^\bullet) \to \widetilde
\H^q(A,M)$ induced by the embedding $\d \ct^\bullet \subset
\ct^\bullet$ are isomorphisms.
\end{remark}

%%%%%%%%%%%%%%%%%%%%%%%%%%%%%%%%%%%%%%%%%%%%%%%%%%%%%%%%%%%%%%%
\section{Extensions and deformations}\label{s4}
%%%%%%%%%%%%%%%%%%%%%%%%%%%%%%%%%%%%%%%%%%%%%%%%%%%%%%%%%%%%%%%

Our \coh\ theory describes extensions and deformations, just as any
\coh\ theory.

\begin{thm}
\label{exts}

\begin{enumerate}

\item
$\widetilde \H^0(A,M) = M^A = \{ m \in M \; | \; a_\lambda m = 0
\;\, \forall a\in A\}$.

\item
The isomorphism classes of extensions
\[
0 \to M \to E \to \nc \to 0
\]
of the trivial $A$-module $\nc$ $(\d$ and $A$ act by zero$)$ by a
conformal $A$-module $M$ correspond bijectively to $\H^0 (A, M)$.

\item
The isomorphism classes of $\nc[\d]$-split extensions
\[
0 \to M \to E \to N \to 0
\]
of conformal modules over a conformal algebra $A$ correspond
bijectively to $\H^1(A, \Chom (N,M))$, where $M$ and $N$ are assumed
to be finite and $\Chom(N,M)$ is the $A$-module of conformal linear
maps from $N$ to $M$. If, in particular, $N =
\nc$ is the trivial module, then there exist no non-trivial
$\nc[\d]$-split extensions.

\item
Let $C$ be a conformal $A$-module, considered as a conformal algebra
with respect to the zero $\lambda$-bracket. Then the equivalence
classes of $\nc[\d]$-split ``abelian'' extensions
\[
0 \to C \to \wti A \to A \to 0
\]
of the conformal algebra $A$ correspond bijectively to
$\H^2(A,C)$.

\item 
The equivalence classes of first-order deformations of a conformal
algebra $A$ $($leaving the $\nc[\d]$-action intact$)$ correspond
bijectively to $\H^2(A,A)$.

\end{enumerate}

\end{thm}

\begin{proof}
1. The computation of $\wti \H^0 (A,M)$ follows directly from the
definition: for $m \in M = \ct^0(A,M)$ and $a \in A$, $(dm)_\lambda
(a) = a_\lambda m$.

2. Given an extension
\[
0 \to M \to E \to \nc \to 0
\]
of modules over a conformal algebra $A$, pick a splitting of this short
exact sequence over $\nc$, \emph{i.e.}, assume that as a complex
vector space, $E \simeq M \oplus \nc = \{(m,n) \; | \; m \in M, n \in
\nc\}$. Define $f \in M$ by writing down the action of $\d$ on
the pair $(m,1) \in E$:
\begin{equation}
\label{zero-cocycle}
\d (m,1) = (\d m + f, 0).
\end{equation}
We claim that $f \in M = \ct^0(A, M)$ defines a zero-cocycle in the
reduced complex $C^\bullet (A, M)$ and thereby a class in $\H^0(A,
M)$.

To see that, define a one-cochain $\gamma \in \ct^1(A, M)$ using the
action of $A$ on $E$:
\begin{equation}
\label{one-cochain.1}
a_\lambda (m,1) = (a_\lambda m + \gamma_\lambda(a), 0)
\end{equation}
for $a \in A$. The conformal antilinearity of $\gamma$:
$\gamma_\lambda (\d a) = -\lambda \gamma_\lambda (a)$, follows from
the fact that $(\d a)_\lambda (m,1) = -\lambda (a_\lambda(m,1))$. The
property $a_\lambda (\d (m,1)) = (\lambda + \d) (a_\lambda(m,1))$ of
the action of $A$ on $E$ expands as
\begin{equation}
\label{df=dg.1}
(df)_\lambda = (\d \gamma)_\lambda,
\end{equation}
which means that $df = 0$ in the reduced complex.

If we choose another splitting $(m,n)'$ of the extension $E$, it will
differ by an element $g \in M$:
\[
(m,1)' = (m+g,1),
\]
so that the new zero-cocycle becomes $f' = f + \d g$, therefore
defining the same cochain in the reduced complex.

If we have two isomorphic extensions and choose a compatible splitting
over $\nc$, we will get exactly the same zero-cocycles $f$
corresponding to them. This proves that isomorphism classes of
extensions give rise to elements of $\H^0(A, M)$.

Conversely, given a cocycle in $C^0(A,M)$, we can choose a
representative $f \in M$ of it to alter the natural $\nc[\d]$-module
structure on $M \oplus \nc$ by adding $f$ to the action of $\d$ on $M
\oplus \nc$ as in \eqref{zero-cocycle}. This will obviously extend to an
action of the free commutative algebra $\nc[\d]$. We can also alter
the natural $A$-module structure by adding $\gamma$ to the action of
$a \in A$ as in \eqref{one-cochain.1}, where $\gamma$ is a solution of
Equation~\eqref{df=dg.1}, which means that $f$ is a cocycle in the
reduced complex. This action will be conformally linear in $(m,n)$,
because of \eqref{df=dg.1}, and antilinear in $A$, because of the
conformal antilinearity of $\gamma$. This action will define an
$A$-module structure on $M \oplus \nc$, because $d \gamma = 0$, which
follows from \eqref{df=dg.1} and the fact that $\nc[\d]$ acts freely on
basic two-cochains.

By construction the natural mappings $M \to M\oplus \nc$ and $M
\oplus \nc \to \nc$ will be morphisms of $\nc[\d]$- and $A$-modules.

This construction of a new conformal module structure on $M \oplus
\nc$ involved a number of choices. The choice of a different
representative $f' = f + \d g$ defines an isomorphism of the two
$\nc[\d]$-module structures on $M \oplus \nc$, which automatically
becomes an isomorphism of the corresponding $A$-module structures,
because the corresponding $\gamma$'s are unique. The one-cochain
$\gamma$ is uniquely determined by $f$, because $\nc[\d]$ acts freely
on the space $\ct^1 (A, M)$ of basic one-cochains.

%\begin{remark}
%This proof does not generalize to arbitrary extensions $0 \to M \to E
%\to N \to 0$, unless $N$ is finite-dimensional over $\nc$, because the
%module $\Hom(N,M)$ is not conformal, in general. The problem is that
%in this case, $\nc[\d]$ will not act freely on the space of nonlocal
%one-cochains, that is, those which are formal power series in
%$\lambda$'s.
%\end{remark}

3. We will adjust the proof of Part 2 to the new situation. Given a
$\nc[\d]$-split extension
\[
0 \to M \to E \to N \to 0
\]
of modules over a conformal algebra $A$, pick a splitting of the short
exact sequence over $\nc[\d]$, \emph{i.e.}, assume that as a
$\nc[\d]$-module, $E \simeq M \oplus N = \{(m,n) \; | \; m \in M, n
\in N\}$. We are going to construct a reduced one-cochain with coefficients in
$\Chom(N,M)$ out of this data. Note that such cochains are linear maps
$\gamma = \gamma_\la (a)_\mu$ from $A \otimes N$ to $M$ depending on
two variables $\la$ and $\mu$, considered modulo $\la - \mu$. Note that
$\gamma_\la(a)_\mu \mod (\la - \mu)$ is fully determined by the
restriction $\gamma_\la (a)_\la$ to the diagonal $\la = \mu$. Define a
one-cochain $\gamma \in C^1(A, \Chom(N,M))$ using the action of $A$ on
$E$:
\begin{equation}
\label{one-cochain}
a_\la (m,n) = (a_\la m + \gamma_\la (a)_\la n, a_\la n)
\end{equation}
for $a \in A$. The conformal antilinearity of $\gamma$:
$\gamma_\lambda (\d a)_\la = -\lambda \gamma_\lambda (a)_\la$, follows
from the fact that $(\d a)_\lambda (m,n) = -\lambda
(a_\lambda(m,n))$. The property $a_\lambda (\d (m,n)) = (\lambda + \d)
(a_\lambda(m,n))$ of the action of $A$ on $E$ expands as
\begin{equation}
\label{df=dg}
(\d \gamma)_\lambda = 0,
\end{equation}
which means that $\gamma_\la (a)_\la$ is a conformal linear map $N \to
M$. Finally, the module property \eqref{module} for elements in $E$
implies that $d \gamma = 0$.

\begin{sloppypar}
If we choose another $\nc[\d]$-splitting $(m,n)'$ of the extension
$E$, it will differ by an element $\beta \in \Hom_{\nc[\d]} (N,M)$:
\[
(m,n)' = (m+\beta(n),n).
\]
$\Hom_{\nc[\d]} (N,M)$ may be identified with the degree-zero part of
$\Chom (N,M)$, so that the new one-cocycle becomes $\gamma' =
\gamma + d \beta$, therefore defining the same cohomology
class.
\end{sloppypar}

If we have two isomorphic extensions and choose a compatible splitting
over $\nc[\d]$, we will have exactly the same one-cocycles $\gamma$
corresponding to them. This proves that isomorphism classes of
extensions give rise to elements of $\H^1(A, \Chom (N,M))$.

Conversely, given a cohomology class in $\H^1(A,\Chom(N,M))$, we can
choose a representative $\gamma \in C^1(A,\Chom(N,M))$ of it to alter
the natural $A$-module structure on $M \oplus N$ by adding $\gamma$ to
the action of $A$ on $M \oplus N$ as in
\eqref{one-cochain}. This action will be conformally linear in $(m,n)$,
because of \eqref{df=dg}, and antilinear in $A$, because of the
conformal antilinearity of $\gamma$. This action will define an
$A$-module structure on $M \oplus N$, because $d \gamma = 0$ after the
restriction to $\mu = \la_1 + \la_2$ in $\ct^2 (A, \Chom(N,M))$.

By construction the natural mappings $M \to M\oplus N$ and $M
\oplus N \to N$ will be morphisms of $\nc[\d]$- and $A$-modules.

This construction of a new conformal module structure on $M \oplus N$
is independent on the choice of a different representative $\gamma' =
\gamma + d \beta$, because it defines an isomorphic structure of an
$A$-module on $M \oplus N$.

Finally, if $N = \nc$, then $\Chom(\nc,M) = 0$, and therefore, there
are no split extensions.

4. Given a $\nc[\d]$-split extension of a conformal algebra $A$ by a
module $C$, choose a splitting $\widetilde A = C \oplus A$ thereof.
Then the bracket in $\widetilde A$
\[
[(0,a)_\lambda (0,b)] = (c_\lambda (a,b), a_\lambda b) \qquad \text{for
$a,b \in A$}
\]
defines a sesquilinear map $c\colon A \otimes A \to C[\lambda]$, which
we may combine with the natural mapping
\begin{align*}
C[\lambda] & \to C[\lambda_1,\lambda_2] /(\d + \la_1 +\la_2), \\
p(\lambda) & \mapsto p(\lambda_1),
\end{align*}
to get the composite mapping, denoted $c_{\lambda_1,\lambda_2}$. It
defines a two-cochain, because it is obviously skew and $(c_\lambda
(\partial a, b), -\la a_\lambda b) = [(0,\partial a)_\lambda (0,b)] =
[\partial (0,a)_\lambda (0,b)] = - \lambda [(0,a)_\lambda (0,b)] =
(-\lambda c_\lambda (a,b), -\lambda a_\lambda b)$, which implies
$c_{\lambda_1,\lambda_2} (\partial a,b) = -\lambda_1
c_{\lambda_1,\lambda_2} (a,b)$, and similarly,
$c_{\lambda_1,\lambda_2} (a,\partial b) \linebreak[1] = \linebreak[0]
-\lambda_2 c_{\lambda_1,\lambda_2} (a,\linebreak[0] b) \mod (\partial
+ \lambda_1 + \lambda_2)$. In fact, this two-cochain $c$ is a cocycle:
\begin{multline*}
dc = a_{\lambda_1} c_{\lambda_2,\lambda_3}(b,c) - b_{\lambda_2}
c_{\lambda_1,\lambda_3}(a,c) + c_{\lambda_3} c_{\lambda_1,\lambda_2}(a,b)
\\
- c_{\lambda_1+\lambda_2,\lambda_3} (a_{\lambda_1} b, c)
+ c_{\lambda_1+\lambda_3,\lambda_2} (a_{\lambda_1} c, b)
- c_{\lambda_2+\lambda_3,\lambda_1} (b_{\lambda_2} c, a) = 0.
\end{multline*}
This is just because the corresponding three-term relation, the Jacobi
relation, is satisfied in $\widetilde A$. 

The construction of $c$ assumed the choice of a splitting $\widetilde
A = C \oplus A$. A different splitting would differ by a mapping
$f\colon A \to C$, which can be thought of as $f\colon A \to
C[\lambda] / (\d + \la)$, which would contribute by $df$ to $c$.

Thus, any extension determines a \coh\ class in $\H^2 (A,
\linebreak[0] C)$.  The above arguments can be traced back to show that a
class in the \coh\ group defines an extension.

5. Let $D = \nc[\epsilon]/(\epsilon^2)$ be the algebra of \emph{dual
numbers}. Then a \emph{first-order deformation} of a conformal algebra
$A$ is the structure of a conformal algebra over $D$ on $A \otimes D$,
so that the map $A \otimes D \to A$, $a \otimes p(\epsilon) \mapsto
p(0) \cdot a$, is a morphism of conformal algebras and the action of
$\d$ on $A \otimes D$ is induced from that on the first factor. This
means classes of first-order deformations are in bijection with
classes of $\nc[\d]$-split abelian extensions of $A$ with the
$A$-module $A$ in the sense of Part 2 of this theorem. Therefore, they
are classified by $\H^2(A,A)$.
\end{proof}

%%%%%%%%%%%%%%%%%%%%%%%%%%%%%%%%%%%%%%%%%%%%%%%%%%%%%%%%%%%%%%%%%%%%%%%%
\section{Homology}\label{s5}
%%%%%%%%%%%%%%%%%%%%%%%%%%%%%%%%%%%%%%%%%%%%%%%%%%%%%%%%%%%%%%%%%%%%%%%

%\subsection{Chains}

%\begin{sloppypar}
Dualizing the \coh\ theory we have defined above, the space
$\ct_n(A,M)$ of
\emph{$n$-chains of a conformal algebra $A$ with coefficients in a
conformal module $M$ over it} is defined as the quotient of
\[
A^{\otimes n} \otimes \Hom( \nc[\lambda_1, \dots, \lambda_n], M),
\]
where $\Hom (\nc[\lambda_1, \dots, \lambda_n], M)$ is the space of
$\nc$-linear maps from the space of polynomials to the module $M$, by
the following relations:
\begin{enumerate}
\item $a_1 \otimes \dots \otimes \partial a_i \otimes \dots \otimes a_n
\otimes \phi 
= - a_1 \otimes \dots \otimes a_i \otimes \dots \otimes a_n
\otimes T_i \phi$, where $(T_i \phi) (f) = \phi (\lambda_i f)$;
\item $a_1 \otimes \dots \otimes a_i \otimes \dots \otimes a_j \otimes 
\dots \otimes a_n \otimes \phi = - a_1 \otimes \dots \otimes a_j \otimes \dots \otimes a_i \otimes 
\dots \otimes a_n \otimes \tau_{ij}^*\phi$, where
$(\tau_{ij}^* \phi) (f(\lambda_1, \dots, \lambda_i, \dots,
\lambda_j, \dots, \lambda_n)) = \phi (f(\lambda_1, \dots, \lambda_j, \dots,
\lambda_i, \dots, \lambda_n))$.
\end{enumerate}
One can also define a differential which takes $n$-chains to
$(n-1)$-chains as follows:
\begin{multline*}
\delta (a_1 \otimes \dots \otimes a_n \otimes \phi )
\\
\begin{split}
& = \sum_{i=1}^{n} (-1)^{i+1} p_i (a_1 \otimes \dots \otimes \widehat
a_i
\otimes \dots \otimes a_{n} \otimes {a_i}_{\lambda_i} \phi )
\\
&
%\begin{split}
+ \sum^{n}_{\substack{i,j = 1\\i < j}}
(-1)^{i+j} p_{ij} ([{a_i}_{\lambda_i} a_j] \otimes a_1 \otimes \dots \otimes 
\widehat a_i \otimes\dots \otimes \widehat a_j \otimes\dots \otimes a_{n}
\otimes \phi),
%\end{split}
\end{split}
\end{multline*}
where $p_i$ is the natural pairing map $\nc[\lambda_i] \otimes
\Hom (\nc[\lambda_1, \dots, \lambda_n], M) \linebreak[1] \to \linebreak[0]
\Hom (\nc[\lambda_1, \dots, \linebreak[0] \widehat
\lambda_i, \dots, \lambda _n], M)$ and $p_{ij}$ is the pairing
$\nc[\lambda_i] \linebreak[2] \otimes \linebreak[1] \Hom (
\nc[\lambda_1, \linebreak[0] \dots, \linebreak[1] \lambda_n], M)
\linebreak[1] \to \linebreak[0] \Hom( \nc[\lambda_i+\lambda_j,
\lambda_1,
\dots, \widehat
\lambda_i, \dots, \widehat \lambda_j, \dots , \lambda _n], M)$.
Similar computations to those in the cochain case show that the
operator $\delta$ is well-defined and $\delta^2=0$.
%\end{sloppypar}

One can define \emph{basic homology} $\wti \H_\bullet (A,M)$ as the
homology of the chain complex and \emph{reduced homology} as the
homology of the subcomplex $C_\bullet (A,M)$ of $\d$-invariant
chains, where $\d$ acts as
\[
\d (a_1 \otimes \dots \otimes a_n \otimes \phi ) 
= a_1 \otimes \dots \otimes a_n \otimes (\d \phi -\sum_{i=1}^n T_i
\phi),
\]
where $(\d \phi)(f) = \d (\phi (f))$, $f \in \nc[\la_1, \dots \la_n]$.
There are obviously natural pairings $\wti \H_q (A,M^*) \otimes \wti
\H^q (A,M) \to \nc$ and $\H_q (A,M^*)
\otimes \H^q (A,M) \to \nc$ for $q \ge 0$, where $M^* = \Hom_\nc (M,
\nc)$ is the linear dual space with a natural structure of an
$A$-module:
\begin{align*}
(\d f) (m) & = - f(\d m),\\
(a_\la f) (m) & = -f(a_\la m)
\end{align*}
for $f \in M^*$, $m \in M$, and $a \in A$. One expects these pairings
to be perfect, when, for instance, either of the (co)homology spaces
is finite-dimensional.

%%%%%%%%%%%%%%%%%%%%%%%%%%%%%%%%%%%%%%%%%%%%%%%%%%%%%%%%%%%%%%%%%%%%%%%
\section{Exterior multiplication, contraction, and module structure}
\label{s7}
%%%%%%%%%%%%%%%%%%%%%%%%%%%%%%%%%%%%%%%%%%%%%%%%%%%%%%%%%%%%%%%%%%%%%%%

For any $u \in \wti C^m (A,\nc)$, where $\nc$ is the
one-dimensional space with the zero action of $A$, let $\epsilon (u)$
be the operator of \emph{exterior multiplication} on 
$\wti C^\bullet(A,M)$:
\begin{multline*}
(\epsilon(u) \gamma)_{\lambda_1, \dots, \lambda_{m+n}} (a_1, \dots,
a_{m+n})
\\
= \sum_{\pi \in S_{m+n}} \frac{\sign \pi}{m!\; n!}\;
u_{\lambda_{\pi(1)},
\dots, \lambda_{\pi(m)}} (a_{\pi(1)},
\dots, a_{\pi(m)})
\\
\times 
\gamma_{\lambda_{\pi(m+1)}, \dots, \lambda_{\pi(m+n)}}
(a_{\pi(m+1)}, \dots, a_{\pi(m+n)}).
\end{multline*}
Define also a wedge product $u \wedge \gamma = \epsilon (u)
\gamma$ on $\ct^\bullet(A,\nc)$. It is clear that $\epsilon (u \wedge v) =
\epsilon (u) \epsilon (v)$ for any $u,v \in V$, therefore, we have a
graded commutative associative algebra structure on $\wti C^\bullet
(A, \nc)$, along with a $\wti C^\bullet (A, \nc)$-module structure on
$\wti C^\bullet (A, M)$.

%\begin{sloppypar}
Similarly, for any chain $v = a_1 \otimes \dots \otimes a_n \otimes
\phi \in \ct_n (A, \nc)$, let $\iota (v)$ be the following
\emph{contraction} operator 
$\wti C^m(A,M) \to \wti C^{m-n}(A,M)$, for $m \ge n$:
\[
%\begin{multline*}
(\iota (v) \gamma)_{\lambda_{n+1}, \dots, \lambda_{m}} (a_{n+1},
\dots, a_{m})
%\\
%= \sum_{\pi \in S_m} p \left(\frac{\sign \pi}{(m-n)! \;
%n!} \;
%\pi^* \gamma \otimes \gamma_{\lambda_{\pi(1)}, \dots, \lambda_{\pi(m)}}
%(a_{\pi(1)} \otimes \dots \otimes a_{\pi(m)})
%\right),
= p\bigl(\phi\tt\ga_{\la_1,\dots,\la_m}(a_1,\dots,a_m)\bigr),
%\end{multline*}
\]
where $p$ is the natural pairing 
$\nc[\lambda_1, \dots, \lambda_n]^* 
\otimes 
\nc[\lambda_1, \dots, \lambda_m]
\linebreak[1] \to \linebreak[0]
\nc[\lambda_{n+1}, \linebreak[0] \dots,
\lambda_{m}]$. 
%and $(\pi^* \gamma) (f(\lambda_1, \dots, \lambda_m)) 
%= \gamma (f(\lambda_{\pi(1)}, \dots, \lambda_{\pi(m)}))$.  
Note that for any $u \in \wti C^1(A,\nc)$ and $v \in \ct_1(A,\nc)$,
\[
\epsilon (u) \iota (v) + \iota (v) \epsilon (u) = \iota(v) u.
\]
%\end{sloppypar}

Furthermore for any $a \in A$, define the following structure of a
\emph{module over the conformal algebra $A$ on $\wti C^\bullet (A,M)$}:
\begin{multline*}
(\theta_\lambda(a) \gamma)_{\lambda_{1}, \dots, \lambda_{n}} (a_{1},
\dots, a_{n})
\\
= a_\lambda \gamma_{\lambda_{1}, \dots, \lambda_{n}} (a_{1},
\dots, a_{n})
%\\
- \sum_{i=1}^{n} \gamma_{\lambda_{1}, \dots,
\lambda + \lambda_i, \dots, \lambda_{n}}(a_{1},
\dots, [a_\lambda a_{i}], \dots, a_{n}).
\end{multline*}
Define $\iota_\lambda(a)$ in a similar fashion:
\begin{equation*}
(\iota_\lambda(a) \gamma)_{\lambda_{1}, \dots, \lambda_{n-1}} (a_{1},
\dots, a_{n-1})
= \gamma_{\lambda,\lambda_{1}, \dots, \lambda_{n-1}}(a, a_{1}, \dots, a_{n-1}).
\end{equation*}
Note that every $a \in A$ defines naturally a one-chain
$a\tt\gamma_{\lambda_0} \in \ct_1(A,\nc)$ depending on 
a parameter $\lambda_0$, where
$\gamma_{\lambda_0} (f(\lambda)) = f(\lambda_0)$. Then we have
$\iota_\lambda(a) = \iota (a\tt\ga_\la)$.  The
fundamental identity
\[
d \iota_\lambda + \iota_\lambda d = \theta_\lambda
\]
of classical Lie theory is also valid in the context of conformal
algebras. It also implies $d \theta_\lambda = \theta_\lambda d$.
As in the Lie algebra case, the induced action of $A$
on $\wti\H^\bullet(A,M)$ is trivial, cf.\ Remark~\ref{rextetc}.

%%%%%%%%%%%%%%%%%%%%%%%%%%%%%%%%%%%%%%%%%%%%%%%%%%%%%%
\section{Cohomology of conformal algebras and their annihilation 
Lie algebras}\label{s8}
%%%%%%%%%%%%%%%%%%%%%%%%%%%%%%%%%%%%%%%%%%%%%%%%%%%%%%
\subsection{Cohomology of the basic complex}
%%%%%%%%%%%%%%%%%%%%%%%

Let $A$ be a conformal algebra and $M$ a conformal module over it.
Then $M$ is a module over the annihilation Lie algebra $\g_-=(\Lie
A)_-$, see Section~\ref{s1}.  Let $C^\bullet(\g_-,M)$ be the
Chevalley--Eilenberg complex defining the cohomology of $\g_-$ with
coefficients in $M$. Recall that, by definition (see, \emph{e.g.}, \cite{F}),
$C^n(\g_-,M)$ is the space of skew-symmetric linear functionals
$\ga\colon (\g_-)^{\tt n} \linebreak[0] \to M$ which are {\em continuous}, 
\emph{i.e.},
\begin{equation*}
\ga({a_1}_{m_1} \tt\dotsm\tt {a_n}_{m_n}) = 0
\end{equation*}
for all but a finite number of $m_1,\dots,m_n\in\ZZ_+$, where $a_1,
\dots, a_n \in A$, and ${a_i}_{m_i} \in \g_- = (\Lie A)_- =
A[t]/(\d + \d_t)A[t]$ is
the image of the element $a_i t^{m_i}$.

$C^\bullet(\g_-,M)$ has the following structure of a $\Cset[\d]$-module:
\begin{multline}\label{dgam}
(\d\ga) (a_1 \otimes \dots \otimes a_n) \\
= \d\bigl(\ga (a_1 \otimes \dots\otimes a_n)\bigr) 
- \sum_{i=1}^n \ga (a_1 \otimes \dots \otimes \d a_i \otimes
\dots \otimes a_n),\\
\qquad \ga\in C^n(\g_-,M).
\end{multline}
%%%%%%%%%%%%%%%%%%%%%%%%%%%
\begin{sloppypar}
\begin{thm}\label{ec1}
There is a canonical isomorphism of
complexes $\wti C^\bullet(A,M)$ and $C^\bullet(\g_-,M)$,
compatible with the action of $\Cset[\d]$.
Consequently, the complex $C^\bullet(A,M)$ 
is isomorphic to $C^\bullet(\g_-,M) / \d C^\bullet(\g_-,M)$.
\end{thm}
\end{sloppypar}
\begin{proof}
For a cochain $\al\in \wti C^n(A,M)$, we write
\begin{equation*}%\label{alm1}
\al_{\la_1,\dots,\la_n}(a_1,\dots,a_n)
= \sum_{m_1,\dots,m_n\in\Zset_+} \la_1^{(m_1)} \dotsm \la_n^{(m_n)}
\al_{(m_1,\dots,m_n)}(a_1,\dots,a_n).
\end{equation*}
In terms of the linear maps
\begin{align*}
%\label{alm2}
&\al_{(m_1,\dots,m_n)} \colon A^{\tt n}\to M,
\\
%\label{alm3}
&a_1\tt\dotsm\tt a_n \mapsto \al_{(m_1,\dots,m_n)}(a_1,\dots,a_n),
\end{align*}
the definition of $\wti C^\bullet(A,M)$ translates as follows.
\begin{enumerate}
\item For any $a_1,\dots,a_n\in A$, 
$\al_{(m_1,\dots,m_n)}(a_1,\dots,a_n)$ 
is non-zero for only a finite number of $(m_1,\dots,m_n)$.

\item  $\al_{(m_1,\dots, m_i, \dots, m_n)}(a_1,\dots, \d a_i, \dots,a_n)$
\\  $= -m_i \al_{(m_1,\dots, m_i-1, \dots, m_n)}(a_1,\dots, a_i, \dots,a_n)$.

\item  $\al$ is skew-symmetric with respect to
simultaneous permutations of $a_i$'s and $m_i$'s.
\end{enumerate}
The differential is given by:
\begin{multline*}
(d\gamma)_{(m_1, \dots, m_{n+1})} (a_1,\dots,a_{n+1})
\\
\begin{split}
& = \sum_{i=1}^{n+1}  
(-1)^{i+1} {a_i}_{(m_i)}
\gamma_{(m_1, \dots, \widehat m_i,
\dots, m_{n+1})} (a_1,\dots, \widehat a_i,\dots, a_{n+1})
\\
&
\begin{split}
+ \sum^{n+1}_{\substack{i,j = 1\\i < j}} \sum_{k=0}^{m_i}
(-1)^{i+j} \binom{m_i}{k}
\gamma_{(m_i+m_j-k, m_1,\dots,\widehat m_i,\dots,\widehat m_j,\dots,m_{n+1})} 
({a_i}_{(k)} a_j, a_1, &
\\
\dots, \widehat a_i,\dots, 
\widehat a_j, \dots, a_{n+1})&.
\end{split}
\end{split}
\end{multline*}
Define linear maps $\phi^n\colon \wti C^n(A,M)\to C^n(\g_-,M)$ 
by the formula
\begin{equation*}
(\phi^n\al)({a_1}_{m_1} \tt\dotsm\tt {a_n}_{m_n})
= \al_{(m_1,\dots,m_n)}(a_1,\dots,a_n).
\end{equation*}
They are well-defined due to above condition 2.  Clearly, $\phi^n$ are
bijective and, using \eqref{eq:1.3}, it is easy to see that
$\phi^{n+1}\circ d = d\circ\phi^n$.  Moreover, $\phi^n\circ\d =
\d\circ\phi^n$, where $\d$ acts on $\wti C^\bullet(A,M)$ via
\eqref{dwtic} and on $C^\bullet(\g_-,M)$ via \eqref{dgam}.
\end{proof}
%%%%%%%%%%%%%%%%%%%%%%%%%%%
\begin{crl}\label{cec1}
$\wti \H^\bullet(A,M) \simeq \H^\bullet(\g_-, M)$.
\end{crl}
%%%%%%%%%%%%%%%%%%%%%%%%%%%
\begin{remark}\label{rhomol}
Similar results hold for homology. To a chain
$a_1\tt\dotsm\tt a_n \tt \phi\in \ct_n(A,M)$
($a_i\in A$, $\phi\in \Hom (\Cset[\la_1,\dots,\la_n], M)$)
we associate the chain
\begin{equation*}
\langle \phi, {a_1}_{\la_1}\tt\dotsm\tt {a_n}_{\la_n} \rangle 
 \in C_n(\g_-,M).
\end{equation*}
In other words, 
$a_1\tt\dotsm\tt a_n \tt
\bigl(\d_{\la_1}^{(m_1)} \dotsm \d_{\la_n}^{(m_n)}
|_{\la_1=\dots=\la_n=0}\bigr)$
corresponds to 
${a_1}_{m_1}\tt\dotsm\tt {a_n}_{m_n}$.
\end{remark}
%%%%%%%%%%%%%%%%%%%%%%%%%%%
\begin{remark}\label{rextetc}
One can easily see that the exterior multiplication, contraction,
module structure, etc., of Section \ref{s7} are equivalent to the
corresponding notions for the annihilation Lie algebra $\g_-$.
For example, if $\th(a_m)$ denotes the action of $a_m\in\g_-$
on $C^\bullet(\g_-,M) \simeq \ct^\bullet(A,M)$, then 
\[
\th_\la(a) = \sum_{m\in\Zset_+} \la^{(m)} \th(a_m).
\]
In particular, the action of $A$ on $\wti\H^\bullet(A,M)$ is trivial.
\end{remark}
%%%%%%%%%%%%%%%%%%%%%%%%%%%
\subsection{Cohomology of the reduced complex}

Now we assume that $M$ is a free $\Cset[\d]$-module: 
$M=\Cset[\d]\otimes_\Cset U$
for some vector space $U$. Then the $\g_-$-module 
$V_-=V(M)_-$ is just $U[t]$ with 
\[
a_{m}(u t^n) = \sum_{j=0}^m \binom{m}{j}\, (a_{(j)}u) \, t^{m+n-j}, \quad
\d (u t^n) = -n ut^{n-1},
\]
for $u\in U$, $a\in A$, see Section~\ref{s1}.
In terms of the usual generating series 
$a_\la = \sum_{m\ge0} \la^{(m)} a_{m}$, this can be rewritten as
\[
a_\la (u t^n) = (a_\la u)\, t^n e^{t \la}.
\]
%Recall also that $V_-$ is an $A$-module\footnote{Of course this module
%is not conformal.} where $\Cset[\d]$ acts by $\d=-\d_t$.

%%%%%%%%%%%%%%%%%%%%%%%%%%%
\begin{thm}\label{ec2}
\begin{sloppypar}
If $A$ is a conformal algebra and $M$ a conformal module which is free
as a $\Cset[\d]$-module, then the complex $C^\bullet(A,M)$ is
isomorphic to the subcomplex $C^\bullet(\g_-,V_-)^\d$ of
$\d$-invariant cochains in $C^\bullet(\g_-, V_-)$.
\end{sloppypar}
\end{thm}
%%%%%%%%%%%%%%%%%%%%%%%%%%%
\begin{proof}
Let $\be\in C^n(\g_-,V_-)$.
As in the proof of Theorem \ref{ec1}, consider the generating series
\begin{multline}\label{blat}
\be_{\la_1,\dots,\la_n;t}(a_1,\dots,a_n)
\\
= \sum_{m_1,\dots,m_n\in\Zset_+} \la_1^{(m_1)} \dotsm \la_n^{(m_n)}
\be({a_1}_{m_1} \tt\dotsm\tt {a_n}_{m_n}).
\end{multline}
By Equation~\eqref{dgam}, $\d$ acts on $\be_{\la_1,\dots,\la_n;t}$ as
$-\d_t+\sum \la_i$. Hence $\be$ is $\d$-invariant, iff 
\begin{equation}\label{psin2}
\be_{\la_1,\dots,\la_n;t}(a_1,\dots,a_n)
= \ga_{\la_1,\dots,\la_n}(a_1,\dots,a_n) \, e^{t \sum \la_i}
\end{equation}
where $\ga_{\la_1,\dots,\la_n} = \be_{\la_1,\dots,\la_n;t}|_{t=0}$
takes values in $U$.
Identifying $U$ with $1\otimes U \subset M$, we can consider $\ga$
as an element of $\wti C^n(A,M)$. It is easy to check that 
$\be\mapsto \ov\ga := \ga \mod(\d+\sum\la_i)$ 
is a chain map from $C^\bullet(\g_-,V_-)$
to $C^\bullet(A,M)$.

Conversely, for $\overline\ga\in C^n(A,M)$ choose a representative
$\ga\in \wti C^n(A,M)$ such that $\overline\ga = \ga \mod(\d+\sum\la_i)$.
Define $\be\in C^n(\g_-,V_-)^\d$ by 
(\ref{blat}, \ref{psin2}) with $\d$ substituted by
$-\d_t$ in 
$\ga_{\la_1,\dots,\la_n}(a_1,\dots,a_n) \in M=U[\d]$.
Then clearly, $\be$ is independent of the choice of $\ga$.

The correspondence $\be\leftrightarrow\ov\ga$ establishes an isomorphism
between $C^\bullet(\g_-,V_-)^\d$ and $C^\bullet(A,M)$.
%This completes the proof.
\end{proof}

\begin{remark}
\label{local}
Identifying $C^\bullet(A,M)$ with $C^\bullet(\g_-,M) / \d C^\bullet(\g_-,M)$,
we can rewrite \eqref{psin2} as 
\begin{multline*}\label{psin}
\be({a_1}_{m_1} \tt\dotsm\tt {a_n}_{m_n})
\\
= \sum_{k_1,\dots,k_n\in\Zset_+}
\binom{m_1}{k_1} \dotsm \binom{m_n}{k_n}
\ga({a_1}_{k_1} \tt\dotsm\tt {a_n}_{k_n}) 
\, t^{\sum m_i - \sum k_i},
\end{multline*}
($a_i\in A, m_i\in\Zset_+$).
\end{remark}

\subsection{Cohomology of conformal algebras and formal distribution
Lie algebras}
\label{relation}

Let $\gtg$ be the maximal formal distribution Lie algebra
corresponding to a conformal algebra $A$ (see Section~\ref{s1}).
Suppose $\gamma \in C^n(A,M)$. The following formula defines an
$n$-cochain $\wti \gamma$ on the Lie algebra $\gtg$:
\begin{multline*}
\wti \gamma (a_1 f_1(t), \dots, a_n f_n(t)) \\
= \Res_{\la_1, \dots, \la_n} \gamma_{\d_1, \dots , \d_n} (a_1, \dots,
a_n) \de (\la_1-\la_2) \dots \de (\la_1-\la_n) f_1 (\la_1) \dots f_n
(\la_n),
\end{multline*}
where $a_i \in A$, $f_i \in \nc[t]$, $\d_i = \d / \d \la_i$, and when
substituting $\d$ into a polynomial, one has to use the divided powers
$\d^{(k)} = \d^k /k!$. This formula is equivalent to the one from
Remark~\ref{local}, where $m_i$'s are now allowed to take negative
values. This correspondence defines a morphism of complexes and,
therefore, \coh.

%%%%%%%%%%%%%%%%%%%%%%%%%%%%%%%%%%%%%%%%%%%%%%%%%%%%%
\section{Cohomology of the Virasoro conformal algebra}\label{s9}
%%%%%%%%%%%%%%%%%%%%%%%%%%%%%%%%%%%%%%%%%%%%%%%%%%%%%

The conformal algebra with one free generator $L$ as a
$\nc[\partial]$-module and $\la$-bracket
\[
[L_\lambda L] = (\partial + 2 \lambda) L
\]
is called the \emph{Virasoro conformal algebra} $\vir$, 
cf.\ Example~\ref{ex:1.2}.

%%%%%%%%%%%%%%%%%%%%%
\subsection{Cohomology of $\vir$ with trivial coefficients}
%%%%%%%%%%%%%%%%%%%%%

Here we will
compute the \coh\ of $\vir$ with trivial
coefficients $\nc$, where both $\partial$ and $L$ act by zero. 

\begin{thm}\label{hvirc}
For the Virasoro conformal algebra $\vir$,
\begin{align*}
\dim \wti \H^q(\vir, \nc) &= \begin{cases}
1 & \text{if $q=0$ or $3$},\\
0 & \text{otherwise},
\end{cases}
\\
\intertext{and}
\dim \H^q(\vir, \nc) &= \begin{cases}
1 & \text{if $q=0$, $2$, or $3$},\\
0 & \text{otherwise}.
\end{cases}
\end{align*}
\end{thm}

\begin{proof}
Let us first identify the \coh\ complex. An $n$-cochain $\gamma$ in
this case is determined by its value on $L^{\otimes n}$:
\[
P(\lambda_1, \dots, \lambda_n) =
\gamma_{\lambda_1, \dots, \lambda_n} (L, \dots, L).
\]
Obviously, $P(\lambda_1, \dots, \lambda_n)$ is a skew-symmetric
polynomial with values in $\nc$. The differential is then determined
by the following formula:
\begin{multline*}
(dP)(\lambda_1, \dots, \lambda_{n+1}) 
\\
= \sum^{n+1}_{\substack{i,j = 1\\i < j}} (-1)^{i+j} 
(\lambda_i - \lambda_j) P(\lambda_i+\lambda_j,
\lambda_1, \dots, \widehat\lambda_i, \dots, \widehat\lambda_j, \dots,
\lambda_{n+1}).
\end{multline*}
This describes the complex $\wti C^\bullet$. The complex $C^\bullet$
producing the \coh\ of $\vir$ is nothing but the 
quotient of $\wti C^\bullet$ by
the ideal spanned by $\sum_{i=1}^n \lambda_i$ in each degree $n$. In
other words, $C^n$ is the space of regular (polynomial) functions on
the hyperplane $\sum_{i=1}^n \lambda_i = 0$ in $\nc^n$ which are skew
in the variables $\lambda_1, \dots, \lambda_n$. This complex appeared
as an intermediate step in Gelfand--Fuchs's 1968 computation \cite{GF}
of the \coh\ of the Virasoro Lie algebra, and the \coh\ of $C^\bullet$
was computed therein.

Consider the following homotopy operator $\wti C^q \to \wti C^{q-1}$
\[
k(P) = \left. (-1)^q \frac{\partial P}{\partial
\lambda_q}\right|_{\lambda_q = 0}.
\]
A straightforward computation shows that $(dk + kd) P = (\deg P
\nolinebreak - \nolinebreak q) P$ for $P \in \wti C^q$, where $\deg P$ is the
total degree of $P$ in $\lambda_1, \dots, \lambda_q$. Thus, only those
homogeneous cochains whose degree as a polynomial is equal to their
degree as a cochain contribute to the \coh\ of $\wti C^\bullet$. These
polynomials must be skew and therefore divisible by $\Lambda_q
\linebreak[3] =
\linebreak[2] \prod_{i <j} (\lambda_i \nolinebreak - \nolinebreak
\lambda_j)$, whose polynomial degree is $q(q-1)/2$. The quadratic
inequality $q(q-1)/2 \le q$ has $q=0$, 1, 2, and 3 as the only
integral solutions. For $q=0$, the whole $\wti C^0 = \nc$ contributes
to $\H^0(\wti C^\bullet)$.  For $q=1$, the only polynomial of degree 1
is $\lambda_1$, up to a constant factor. $d \lambda_1 =
\lambda_2^2 - \lambda_1^2$, which is the only skew polynomial of
degree 2 in two variables. This shows that $\wti\H^1 = \wti\H^2 = 0$. Finally,
for $q=3$, the only skew polynomial of degree 3 in 3 variables is
$\Lambda_3$, up to a constant. It is easy to see that this polynomial
represents a non-trivial class in the \coh. Indeed, it is closed,
because a skew-symmetric function in four variables has a degree at
least 6, which is greater than $\deg (d \Lambda_3) = 4$. And
$\Lambda_3$ is not a coboundary, because it can be the coboundary of a
two-cochain of degree 2, which must be a constant factor of
$\lambda_2^2 - \lambda_1^2 = d \lambda_1$, whose coboundary is zero.

The computation of the \coh\ of the quotient complex $C^\bullet$ is
based on the short exact sequence
\begin{equation}
\label{long}
0 \to \d \wti C^\bullet \to \wti C^\bullet \to C^\bullet \to 0.
\end{equation}
By definition, $\d \wti C^0 = 0$. To find the
\coh\ of $\d \wti C^\bullet$, define a homotopy $k_1\colon  \d \wti C^q \to
\d \wti C^{q-1}$ as $k_1(\d P) =
\d k(P)$, where $\d = \sum_i \lambda_i$ and $P
\in \wti C^q$. Then $(dk_1 + k_1 d) \d P = (\deg P - q) \d
P$. As in the previous paragraph, this implies that $\deg P = q = 0$,
1, 2, or 3. Up to constant factors, the only polynomials in $\d
\wti C^\bullet$ with this property are $P_1 = \lambda_1^2$ for $q=1$, $P_2
= (\lambda_1 + \lambda_2)(\lambda_1^2 - \lambda_2^2)$ for $q=2$, and
$P_3 = (\lambda_1 + \lambda_2 + \lambda_3) \Lambda_3$ for $q=3$. One
computes: $dP_1 = - P_2$ and $d P_3 = 0$. Therefore $\H^q(\d
\wti C^\bullet) = 0$ for all $q$ but $q=3$, where it is one-dimensional
with the generator $P_3$.

%\begin{sloppypar}
Thus, the long exact sequence of \coh\ associated with \eqref{long}
looks as follows:
\[
\begin{CD}
0 @>>>  0 @>>> \nc  @>>> \H^0(\vir, \nc) @>>> \\
  @>>>  0  @>>>  0  @>>>  \H^1(\vir, \nc)@>>> \\
  @>>>  0  @>>>  0  @>>>  \H^2(\vir, \nc)@>>>\\
  @>>>  \nc P_3  @>>>  \nc \Lambda_3  @>>>  \H^3(\vir, \nc) @>>>\\
  @>>>  0  @>>>  0  @>>>  \H^4(\vir, \nc) @>>>\\
  @>>>  0 @>>>  \dots
\end{CD}
\]
We see that $\H^0(\vir, \nc) = \nc$ and $\H^q(\vir, \nc) = 0$ for $q =
1, 4, 5, 6, \dots$ and $\H^3(\vir, \nc) = \nc \Lambda_3$ and $\H^2(\vir,
\nc) \linebreak[1] = \linebreak[0] \nc (\lambda_1^3 - \lambda_2^3)$, because $d
(\lambda_1^3 - \lambda_2^3) = P_3$.
%\end{sloppypar}
\end{proof}

\begin{remark}
In fact, this computation shows that the \coh\ of the Virasoro
conformal algebra is the primitive part of the \coh\ ring of the
Virasoro Lie algebra, in addition to $\nc$ in degree 0. The reduction
of the basic \coh\ to the computation of Gelfand and Fuchs \cite{GF}
might be made using Corollary~\ref{cec1}, but we preferred to use a
direct argument in the proof.
\end{remark}

\begin{remark}
\label{c_a}
Instead of the trivial $\vir$-module $\nc$, consider the module
$\nc_a$, which is the one-dimensional vector space $\nc$ on which all
elements of $\vir$ act by zero, and $\d v = a v$ for $v \in \nc_a$, $a
\ne 0$ being a given complex constant. Then Proposition~\ref{prop-dd}
shows that $\H^q(\d \ct^\bullet) \simeq \wti \H^q(\vir,\nc_a)$ for $q
\ge 0$, and the long exact sequence \eqref{long-exact} combined with
the computation of $\wti \H^\bullet (\vir, \nc_a)$, which is obviously
isomorphic to $\wti \H^\bullet (\vir, \nc)$, provided by
Theorem~\ref{hvirc}, shows that $\H^q(\vir, \nc_a) = 0$ for all $q$.
\end{remark}

%%%%%%%%%%%%%%%%%%%%%
\subsection{Cohomology of $\vir$ with coefficients in $M_{\De,\al}$}
%%%%%%%%%%%%%%%%%%%%%
Recall (Example~\ref{ex:1.2}) that $M_{\De,\al}$ $(\De,\al\in\Cset)$ is
the following $\vir$-module
\begin{displaymath}
  M_{\Delta , \alpha} = \CC [ \partial ]v \, , \quad
  L_\la v= ( \partial + \alpha + \Delta \lambda )v \, .
\end{displaymath}
As in the previous subsection, we identify the space of $n$-cochains
$C^n(\vir,M_{\De,\al})$ with the space of all $\Cset$-valued
skew-symmetric polynomials
in $n$ variables: for any $\ga\in C^n(\vir,M_{\De,\al})$, there is a unique
polynomial $P(\la_1,\dots,\la_n)$ such that
\[
\ga_{\la_1,\dots,\la_n}(L,\dots,L) = P(\la_1,\dots,\la_n) v
\mod (\d + \la_1+\dots+\la_n).
\]
Then the differential is given by the formula
\begin{align*}
(dP)&(\la_1,\dots,\la_{n+1}) =
\\
&= \sum_{i=1}^{n+1} (-1)^{i+1} 
\Bigl( \al - \sum_{j=1}^{n+1} \la_j + \De\la_i \Bigr)
P(\la_1,\dots,\widehat\la_i,\dots,\la_{n+1})
\\
&+ \sum^{n+1}_{\substack{i,j = 1\\i < j}} (-1)^{i+j} 
(\la_i-\la_j) P(\la_i+\la_j,\la_1,\dots,\widehat\la_i,
\dots,\widehat\la_j,\dots,\la_{n+1}) .
\end{align*}
Now we interpret this in terms of the Lie algebra $\Vect\CC$ of regular
vector fields on $\CC$, which is the annihilation algebra of $\vir$,
see Section~\ref{s1}. To $\ga\in C^n(\vir,M_{\De,\al})$ we associate
a linear map $\be\colon \bigwedge^n \Vect\CC \to \CC$ by the formula
\[
\sum_{m_1,\dots,m_n\in\Zset_+} \la_1^{(m_1)} \dotsm \la_n^{(m_n)}
\be(L_{(m_1)} \wedge\dots\wedge L_{(m_n)})
= P(\la_1,\dots,\la_n),
\]
where $L_{(m)} = -t^m\d_t$.

Then the differential is
\begin{multline*}
(d\be)(L_{(m_1)} \wedge\dots\wedge L_{(m_{n+1})})
\\*
\begin{split}
&
\begin{split}
= \sum_{i=1}^{n+1} (-1)^{i+1} 
\bigl( \al\de_{m_i,0} + (\De-1)\de_{m_i,1} \bigr) \,
\be \bigl(L_{(m_1)} \wedge\dots\wedge \widehat{ L_{(m_i)} } 
& \wedge \dotsm
\\
\dotsm & \wedge L_{(m_{n+1})} \bigr)
\end{split}
%
%\\
%&+ \sum^{n+1}_{\substack{i,j = 1\\i \ne j}} (-1)^{i} 
%m_j \de_{m_i,0} \,
%\be \bigl(L_{(m_1)} \wedge\dots\wedge \widehat{ L_{(m_i)} } \wedge\dots\wedge
%L_{(m_j-1)} \wedge\dots\wedge L_{(m_{n+1})} \bigr)
%
\\
& \begin{split}
+ \sum^{n+1}_{\substack{i,j = 1\\i < j}} (-1)^{i+j} 
(1-\de_{m_i,0}) (1-\de_{m_j,0}) \, \be \bigl( [L_{(m_i)}, L_{(m_j)}],
L_{(m_1)} & \wedge\dotsm
\\
\dotsm\wedge \widehat{ L_{(m_i)} } \wedge\dotsm\wedge
\widehat{ L_{(m_j)} } \wedge\dots &\wedge L_{(m_{n+1})} \bigr).
\end{split}
\end{split}
\end{multline*}

Let $\Vect_0\CC$ be the subalgebra of $\Vect\CC$ of vector fields that
vanish at the origin. It is spanned by the elements
$L_{(m)} = -t^m\d_t$, $m\ge1$. Let $U_\De$ be a $1$-dimensional
$\Vect_0\CC$-module on which $L_{(m)}$ acts as $0$ for $m\ge2$
and $L_{(1)}$ acts as a multiplication by $\De$.

%%%%%%%%%%%%%%%%%%%%%%%
\begin{thm}\label{hvirmd}
\begin{enumerate}
\item
$\H^\bullet(\vir, M_{\De,\al}) = 0$ if $\al\ne0$.

\item
$\H^q(\vir, M_{\De,0}) \simeq 
\H^q(\Vect_0\CC, U_{\De-1}) \oplus 
\H^{q-1}(\Vect_0\CC, U_{\De-1})$ for any $q$ $(H^q = 0$
for $q < 0$ by definition$)$.

\item
$\dim \H^q(\vir, M_{\De,0}) = \dim \H^q(\Vect\CC, \CC[t,t^{-1}]
(dt)^{1-\De})$. Explicitly:
\[
\dim \H^q(\vir, M_{1 - (3r^2 \pm r)/2, \, 0}) 
= \begin{cases}
2 & \text{for $q=r+1$},\\
1 & \text{for $q=r,r+2$},\\
0 & \text{otherwise},
\end{cases}
\]
and $\H^q(\vir, M_{\De,0}) = 0$ if $\De \ne 1 - (3r^2 \pm r)/2$ for any 
$r\in\ZZ_+$.
\end{enumerate}
\end{thm}
\begin{proof}
\begin{sloppypar}
We have seen that the complex $C^\bullet(\vir, M_{\De,\al})$ is isomorphic
to $\left( \bigwedge^\bullet \Vect\CC \right)^*$ endowed with
the above non-standard differential. 
Let $\pi\colon \left( \bigwedge^q \Vect\CC  \right)^* \to
\left( \bigwedge^q \Vect_0\CC  \right)^*$
be the restriction map. It is easy to see that in fact $\pi$ is a
chain map from $C^q (\vir, M_{\De,\al})$ to $C^q (\Vect_0\CC,
\linebreak[0] U_{\De-1})$, where we identify $U_{\De-1}=\CC$ as a
vector space.  Define another map $\iota\colon \left( \bigwedge^{q-1}
\Vect\CC
\right)^* \to
\left( \bigwedge^q \Vect_0\CC  \right)^*$ by the formula
\begin{multline*}
(\iota\be)(L_{(m_1)} \wedge\dots\wedge L_{(m_q)})
\\*
= \sum_{i=1}^{q} (-1)^{i+1} \de_{m_i,0} \,
\be\bigl(L_{(m_1)} \wedge\dots\wedge \widehat{ L_{(m_i)} } \wedge\dots\wedge
L_{(m_q)} \bigr) .
\end{multline*}
Then $\iota$ is a chain map from 
$C^{q-1}(\Vect_0\CC, U_{\De-1})$ to $C^q (\vir, M_{\De,\al})$.
\end{sloppypar}

We have a short exact sequence of complexes
\[
0 \to C^{q-1}(\Vect_0\CC, U_{\De-1}) \stackrel{\iota}{\to}
C^q (\vir, M_{\De,\al}) \stackrel{\pi}{\to}
C^q (\Vect_0\CC, U_{\De-1}) \to 0 .
\]
A splitting 
$\phi\colon C^q (\Vect_0\CC, U_{\De-1}) \to 
\left( \bigwedge^q \Vect_0\CC  \right)^*$
is given by the formula
\[
(\phi\be)(L_{(m_1)} \wedge\dots\wedge L_{(m_q)})
= \begin{cases}
\be(L_{(m_1)} \wedge\dots\wedge L_{(m_q)})
 & \text{if all $m_i\ge1$},\\
0 & \text{otherwise.}
\end{cases}
\]
One checks that if $d\be=0$ then $d\phi\be=\al\iota\be$.

Hence, the cohomology long exact sequence associated to the above short
exact sequence of complexes looks as follows:
\begin{align*}
&\xrightarrow{\al\id} \H^{q-1}(\Vect_0\CC, U_{\De-1})
&\xrightarrow{\iota}\, & \H^q(\vir, M_{\De,\al})
\xrightarrow{\pi} \H^q(\Vect_0\CC, U_{\De-1})
\\
&\xrightarrow{\al\id} \H^q(\Vect_0\CC, U_{\De-1})
&\xrightarrow{\iota}\, & \dotsm
\end{align*}

This proves Parts 1 and 2.

Part 3 follows from Part 2 and the results of Feigin and Fuchs, see
\cite[$\S 2.3$]{F}.  (Note that our $U_\De$ is exactly their
$E_{-\De}$.)
\end{proof}

%%%%%%%%%%%%%%%%%%%%%%%%%%%%%%%%%%%%%%%%%%%%%%%%%%%%%%
\section{Cohomology of current conformal algebras}
\label{scc}
%%%%%%%%%%%%%%%%%%%%%%%%%%%%%%%%%%%%%%%%%%%%%%%%%%%%%%
\subsection{Cohomology with trivial coefficients}
%%%%%%%%%%%%%%%%%%%%%%%%

Here we will compute the \coh\ of a current conformal algebra
$\Cur \gtg$ with trivial coefficients for a finite-dimensional
semisimple Lie algebra $\gtg$. Recall from Example~\ref{ex:1.1} 
that the current conformal algebra $\Cur\g$ is
$\CC[\d]\tt\g$ with the $\la$-bracket
\[
[a_\la b] = [a,b] \qquad\text{for }\; a,b\in\g.
\]
The basic complex in this case becomes bigraded, the second
grading given by the total degree in $\lambda_i$, which we will call
the $\lambda$-\emph{degree}, of the restriction of the cochain to the
subspace $\gtg$ of generators of $\Cur \gtg$. The differential respects
the $\lambda$-degree, and therefore the complex splits into the direct
sum of its graded subcomplexes. Let $\wti C^{\bullet}_0 \subset \wti
C^{\bullet}$ be the subcomplex of zero $\lambda$-degree. This
subcomplex is obviously isomorphic to the Chevalley--Eilenberg complex
$C^\bullet (\gtg, \nc)$ of the Lie algebra $\gtg$.

\begin{thm}
\label{curr-triv}
\begin{enumerate}
\item
The embedding $C^{\bullet}(\gtg,\nc) \subset \wti C^{\bullet}$ is a
quasi-iso\-mor\-phism, \emph{i.e.}, it induces an isomorphism on
\coh. Therefore,
\[
\wti \H^\bullet (\Cur \gtg, \nc) \simeq \H^\bullet (\gtg, \nc) \simeq
\Bigl( \bigwedge\nolimits^\bullet \g^* \Bigr)^\g.
\]

\item
For $q \ge 0$
\[
\H^q (\Cur \gtg, \nc) \simeq \H^q (\gtg, \nc) \oplus \H^{q+1} (\gtg, \nc).
\]

\end{enumerate}
\end{thm}

\begin{proof}
1. According to Theorem~\ref{ec1}, the complexes $\wti C^\bullet(\Cur
\g,\nc)$ and $C^\bullet(\g[t],\nc)$ are isomorphic, because $\g[t]$ is the
annihilation subalgebra of $\Cur \g$, see Example~\ref{ex:1.1}. 
Moreover, the part of $\lambda$-degree zero maps
isomorphically to the Chevalley--Eilenberg complex $C^\bullet(\g,\nc)$,
which is the subcomplex of $C^\bullet (\g[t],\nc)$ of cochains
vanishing on $t\g[t]$. Thus, $\wti \H^\bullet (\Cur \gtg, \nc)
\simeq \H^\bullet (\gtg[t], \nc)$, which is isomorphic to $\H^\bullet
(\gtg,\nc)$ via the subcomplex of cochains vanishing on $t\g[t]$ by a
result of Feigin \cite{Fe,Fe2}; see a different proof of Feigin's result
in Section~\ref{s8.2}, which covers the case of non-trivial
coefficients as well. The computation of $\H^\bullet(\g,\nc)$ via the
invariants of the dual exterior algebra is standard, see \emph{e.g}.,
\cite{Fe2}.

2. Consider the long exact sequence \eqref{long-exact}. The mapping
$\H^q(\partial \ct^\bullet) \to \H^q (\ct)$ for $q \ge 1$ is zero,
because the \coh\ of $\H^q (\ct)$ is concentrated in $\lambda$-degree
zero (see the first statement of the Theorem) and the
\coh\ of $\H^q(\partial \ct^\bullet)$ is concentrated in
$\lambda$-degree one (see Proposition~\ref{prop-dd}). The same is true
even for $q =0$, because $\d \ct^0 = 0$ and the degree-zero
differential $d\colon  \ct^0 \to \ct^1$ is zero. Thus
\eqref{long-exact} splits into the short exact sequences
\[
0 \to \H^q (\gtg, \nc) \to \H^q (\Cur \gtg, \nc) \to \H^{q+1} (\gtg, \nc)
\to 0
\]
for each $q \ge 0$.
\end{proof}

\begin{sloppypar}
\begin{remark}
The same argument as in Remark~\ref{c_a} shows that $\H^\bullet (\Cur
\gtg, \nc_a) = 0$, where $\nc_a$ is the one-dimensional $\Cur \gtg$-module
on which $\Cur \g$ acts trivially and $\d$ acts by
a multiplication by $a \ne 0$.
\end{remark}
\end{sloppypar}

%%%%%%%%%%%%%%%%%%%%%%%%
\subsection{Cohomology with coefficients in a current module}
\label{s8.2}
%%%%%%%%%%%%%%%%%%%%%%%%
Let $\g$ be a finite-dimensional simple Lie algebra,
and $U$ a $\g$-module.
Recall (Example~\ref{ex:1.1})
that the current module $M_U$ over $\Cur\g$ is defined as 
$M_U = \CC[\d]\tt U$ with 
\[
a_\la u = au \qquad\text{for }\; a\in\g, \, u\in U.
\]
%%%%%%%%%%%%%%%%%%%%%%%%
\begin{prop}\label{hcurgmu}
$\H^\bullet(\Cur\g, M_U) \simeq \H^\bullet(\g[t], U)$
where the Lie algebra $\g[t]$ acts on the $\g$-module $U$
by evaluation at $t=0$.
\end{prop}
This can be deduced from Theorem \ref{ec2} but we will give a more
direct argument.
\begin{proof}
Since $M_U$ is free over $\CC[\d]$, any cochain $\al\in C^n(\Cur\g, M_U)$
has a unique representative $\mod(\d+\la_1+\dots+\la_n)$ independent of $\d$.
Explicitly, there is a unique 
$\be\colon {\g}^{\tt n} \to \CC[\la_1,\dots,\la_m] \tt U$
such that
\begin{multline*}
\al_{\la_1,\dots,\la_n}(a_1,\dots,a_n)\\
= \be_{\la_1,\dots,\la_n}(a_1\tt\dotsm\tt a_n) \mod(\d+\la_1+\dots+\la_n)
\end{multline*}
for $a_1,\dots,a_n \in \g$.
Now writing
\begin{multline*}
\be_{\la_1,\dots,\la_n}(a_1\tt\dotsm\tt a_n)
\\
= \sum_{m_1,\dots,m_n \in\ZZ_+}
\la_1^{(m_1)}\dotsm\la_n^{(m_n)}
\, \be(t^{m_1} a_1 \wedge\dots\wedge t^{m_n} a_n)
\end{multline*}
we can interpret $\be$ as a cochain $\bigwedge^n \g[t] \to U$,
as in the proof of Theorem~\ref{ec1}.
\end{proof}
%%%%%%%%%%%%%%%%%%%%%%%%

To compute $\H^\bullet(\g[t], U)$, we apply the Hochschild--Serre
spectral sequence (see, \emph{e.g.}, \cite[$\S 1.5.1$]{F}) 
for the ideal $t \g[t]$ of $\g[t]$. The $E_2$ term is
\begin{align}
\label{e2}
E_2^{p,q} &\simeq \H^p\bigl(\g, \H^q(t \g[t], U)\bigr)
\simeq \H^p(\g) \tt \H^q(t \g[t], U)^{\g}
\\*
\notag
&\simeq \H^p(\g) \tt \bigl(\H^q(t \g[t]) \tt U\bigr)^{\g}.
\end{align}
We used that $U$ is a trivial $t \g[t]$-module and that 
$\H^p(\g, U) \simeq \H^p(\g) \tt U^{\g}$
for any module $U$ over a simple Lie algebra $\g$.

Of course, $\H^p(\g)$ is well-known (cf.\ Theorem~\ref{curr-triv}), so
we only need $\H^q(t \g[t])$. The latter can be deduced from a famous
result of Kostant \cite{Ko} (generalized to the affine Kac--Moody
case).

First, we need some notation from \cite{K0}.  Fix a triangular
decomposition $\g = \n_- \oplus\h\oplus \n_+$.  Let $W$, $\De$,
$\De_+$, $\De_l$, $\rho$, $\theta$, $h^\vee$ be respectively the Weyl
group, the set of roots, the set of positive roots, the set of long
roots, the half sum of positive roots, the highest root, and the dual
Coxeter number of $\g$.  Let $\what\g = \g[t,t^{-1}] + \CC K + \CC d$
be the affine Kac--Moody algebra associated to $\g$.  The
corresponding objects for $\what\g$ will be hatted.  For example,
$\what\g = \what\n_- \oplus\what\h\oplus \what\n_+$, where $\what\n_\pm
= t^{\pm 1} \g[t^{\pm 1}] + \n_\pm$, $\what\h = \h + \CC K + \CC
d$. Denote by $\de$ and $\La_0$ the elements of $\what\h^*$ that
correspond to $K$ and $d$ via the isomorphism $\what\h^* \simeq
\what\h$ given by the invariant bilinear form $(\cdot|\cdot)$ of
$\what\g$, normalized by $(\theta|\theta)=2$.  Recall that the simple
roots of $\what\g$ are $\what\al_0 = \de - \theta$, $\what\al_i =
\al_i$ $(1\le i\le l := \rank\g)$, where $\al_i$ are the simple
roots of $\g$.  The element $\what\rho\in\what\h^*$ is defined by the
property $\langle \what\rho, \what\al_i \rangle = 1$ $(0\le i\le l)$,
\emph{i.e.}, $\what\rho = \rho + h^\vee \La_0$. We denote by bar the
projection from $\what\h^*$ onto $\h^*$. Also recall that $\what W = W
\ltimes T$, where $T$ is the group of translations $t_\ga$
$(\ga\in\ZZ\De_l)$ such that $\ov{t_\ga(\la)} = \la + \langle\la,
K\rangle \ga$ for $\la\in\h^*$, ($w\in W$ acts on $t_\ga$ by $w t_\ga
w^{-1} = t_{w(\ga)}$).  For $\what w\in\what W$ we denote its length
by $\ell(\what w)$.  Finally, if $\La\in\h^*$ is a dominant weight, we
denote by $V(\La)$ the irreducible $\g$-module with highest weight
$\La$.

Now we can state
%%%%%%%%%%%%%%%%%%%%%%%%
\begin{lm}\label{htgt}
\begin{enumerate}
\item
As a $\g$-module
\[
\H^q(t \g[t]) \simeq 
\bigoplus_{\what w\in\what W^1, \, \ell(\what w)=q} 
V \bigl(\ov{\what w(\what\rho) - \what\rho} \bigr)
\]
where $\what W^1 :=
\{ \what w\in\what W \,|\, \what w^{-1}\De_+ \subset \what\De_+ \}$.

\item
Equivalently,
\[
\H^q(t \g[t]) \simeq 
\bigoplus_{(w,\ga)\in WT^1, \, \ell(t_\ga w)=q} 
V(w(\rho) - \rho + h^\vee \ga)
\]
where $WT^1 := \{ (w,\ga) \in W\ltimes\ZZ\De_l \,|\, (\ga|\al) \ge0 \;
\forall\al\in \De_+ , \;\; (\ga|\al) \nolinebreak > \nolinebreak 0
\linebreak[0]
\forall\al\in \De_+ \cap w\De_- \}$.
\end{enumerate}
\end{lm}
\begin{proof}
Part 1 is a special case of Theorem 5.14 of Kostant \cite{Ko}
(generalized to the affine Kac--Moody case).  His Lie algebra $\g$
will be the affine Kac--Moody algebra $\what\g$.  We take the
parabolic subalgebra $\u = \g[t] + \CC K + \CC d$ of $\what\g$, then
$\n = t \g[t]$, $\g_1 = \g + \CC K + \CC d$.

Part 2 is standard, using that $\what W = W \ltimes T$ (see \cite{K0}).
\end{proof}
%%%%%%%%%%%%%%%%%%%%%%%%
\begin{lm}\label{uniqla}
For any $\what w\in\what W$, we have:
\begin{enumerate}
\item
$\what\rho - \what w(\what\rho) = 
\sum_{ \be\in \what\De_+\cap\what w\what\De_- } \be$.

\item
$\ell(\what w) = \big| \what\De_+\cap\what w\what\De_- \big|$.

\item
$\what\rho - \what w(\what\rho) \in\ZZ\de$,
iff $\what w = 1$.
\end{enumerate}
\end{lm}
\begin{proof}
Parts 1 and 2 are exercises from \cite[Chap. 3]{K0} and left to
the reader.

Suppose $\what\rho - \what w(\what\rho) = n\de$, $n\in\ZZ$. Then by Part 1,
$n\de\in\ZZ_+\what\De_+$. Since $\what w^{-1}(\de)=\de$, 
applying $\what w^{-1}$
to Part 1, we get $n\de\in\ZZ_+\what\De_-$. Hence $n=0$. But then Parts 1 and 2
imply $\ell(\what w) = 0$, \emph{i.e.}, $\what w = 1$.
\end{proof}
%%%%%%%%%%%%%%%%%%%%%%%%

It follows from Part 3 of the lemma that for any $\La\in\h^*$
there is at most one $\what w\in\what W$ such that 
$\La = \ov{\what w(\what\rho) - \what\rho}$. Define $\ell(\La)$ to be the
length of this $\what w$ if it exists, and $+\infty$ otherwise.
Then we can restate Lemma~\ref{htgt} as follows:
\begin{equation}\label{cohtgt}
\H^q(t \g[t]) \simeq
\bigoplus_{\La\in\h^* , \, \ell(\La)=q} V(\La),
\end{equation}
where $V(\La)$ is a finite-dimensional representation of highest
weight $\La$.
%%%%%%%%%%%%%%%%%%%%%%%%
\begin{thm}\label{cohcurnontr}
Let $\g$ be a finite-dimensional simple Lie algebra with a fixed Cartan 
subalgebra $\h$. Let
$U$ be an irreducible $\g$-module.
Then
\begin{equation*}
\H^n(\Cur\g, M_{U}) 
\simeq \H^n(\g[t], U) 
\simeq \H^{n-\ell^*(U)}(\g).
\end{equation*}
Here $\ell^*(U) = +\infty$ whenever $U$ is infinite-dimensional, 
$\ell^*(U) = \ell (\Lambda^*)$ whenever $U=V(\Lambda)$ is a finite-dimensional
irreducible module with a highest weight $\La$, 
$\La^*$ is the highest weight of the contragredient module $V(\La)^*$,
$\ell(\La)$ is as above, and we agree that $\H^n=0$ for $n<0$ $($including
$n=-\infty)$.
\end{thm}
\begin{proof}
The first isomorphism in the theorem is from Proposition~\ref{hcurgmu}.
To compute $\H^\bullet(\g[t], V(\La))$,
we apply the Hochschild--Serre spectral sequence for the Lie algebra
$\g[t]$, its module $U$, and its ideal $t\g[t]$.

If $U = V(\La)$, then $U^* \simeq V(\La^*)$ and 
Equations (\ref{e2}, \ref{cohtgt}) imply
that the $E_2$ term is
\[
E_2^{p,q} =  
\begin{cases}
\H^p(\g) & \text{for }\; q=\ell(\La^*) < +\infty,
\\
0 & \text{otherwise}.
\end{cases}
\]
Hence the spectral sequence degenerates at $E_2$ and
$\H^n(\g[t], V(\La)) \simeq \H^{n-\ell(\La^*)}(\g)$.

If $U$ is infinite-dimensional, then again by (\ref{e2}, \ref{cohtgt}),
we have $E_2^{p,q} = 0$.
\end{proof}
%%%%%%%%%%%%%%%%%%%%%%%%
\begin{crl}\label{cohgttriv}
\textup{\cite{Fe,Fe2}}.
$\H^\bullet(\g[t]) \simeq \H^\bullet(\g)$ where the isomorphism is
induced from evaluation at $t=0$.
\end{crl}
%%%%%%%%%%%%%%%%%%%%%%%%
\begin{crl}\label{h1,2cur}
For any semisimple $\g$-module $U$:
\begin{enumerate}
\item
$\H^1(\Cur\g, M_U) \simeq \Hom_\g(\g,U)$.
Explicitly, the isomorphism is giv\-en by:
\[
\al_\la(a) = \la \, \varphi(a) \mod(\d+\la)
\]
for $a\in\g$, $\varphi\in \Hom_\g(\g,U)$.

\item
$\H^2(\Cur\g, M_U) \simeq \Hom_\g(\bigwedge^2\g / \g, U)$,
provided that $\g\not\simeq\sl_2$.
Explicitly, the isomorphism is given by: 
\[
\al_{\la_1,\la_2}(a_1,a_2) 
= \la_1\la_2 \, \varphi(a_1,a_2) \mod(\d+\la_1+\la_2)
\]
for $a_1,a_2\in\g$, $\varphi\in \Hom_\g(\bigwedge^2\g / \g, U)$.
\end{enumerate}
\end{crl}
\begin{proof}
It is easy to check that the above formulas indeed give cocycles.
In fact, 2.\ gives a cocycle for any 
$\varphi\in \Hom_\g(\bigwedge^2\g, U)$; however, any
$\varphi\in \Hom_\g(\bigwedge^2\g, \g)$ gives a coboundary.
Next, we use Lemma~\ref{htgt} and the fact that 
$\H^n(\Cur\g, M_U) \simeq \Hom_\g(\H^n(t\g[t]),U)$ for $n=1,2$.

1. The only element of $\what W^1$ of length $1$ is the simple reflection
$r_{\what\al_0}$ with respect to the root $\what\al_0$. Then 
$r_{\what\al_0}(\what\rho)-\what\rho = -\what\al_0 = \theta-\delta$. Hence
$\H^1(t\g[t]) \simeq V(\th)\simeq\g$ as a $\g$-module.

2. All elements of $\what W^1$ of length $2$ are of the form 
$r_{\what\al_0} r_{\what\al_i}$, where $i$ is such that 
$\langle \what\al_i, \what\al_0^\vee \rangle \ne0$. Then
$r_{\what\al_0} r_{\what\al_i}(\what\rho)-\what\rho 
= - \what\al_0 - \what\al_i + 
\langle \what\al_i, \what\al_0^\vee \rangle \what\al_0$.

When $\g=\sl_2$ we get $\H^2(t\g[t]) \simeq V(2\al_1)$, see the next
example.  When $\g=\sl_{l+1}$, $l\ge2$, there are two possibilities
for $i$: either $i=1$ or $i=l$; then $\H^2(t\g[t]) \simeq
V(2\th-\al_1) \oplus V(2\th-\al_l)$.  For $\g\ne\sl_{l+1}$ there is a
unique possibility for $i$ and $\H^2(t\g[t]) \simeq V(2\th-\al_i)$.

In all cases, except for $\sl_2$, one can check that 
$\H^2(t\g[t]) \simeq \bigwedge^2\g / \g$.
\end{proof}
%%%%%%%%%%%%%%%%%%%%%%%%
\begin{ex}\label{cohsl2}
Let $V(m)$ be the unique irreducible $\sl_2$-module of dimension $m+1$.
Then $\dim \H^n(\Cur\sl_2, M_{V(m)}) = 1$ for $m = 2n, 2(n-3)$, 
and $=0$ otherwise. 

Let $\{e,f,h\}$ be the standard basis of $\sl_2$. Then the module
$V(2n)$ is isomorphic to $\Sym^n \sl_2 / (h^2-4ef)$. Note that
$\Sym^\bullet \sl_2 / (h^2-4ef)$ is the coordinate ring of the
nilpotent cone of $\sl_2$. This description of $V(2n)$ allows us to
give an explicit formula for the cocycles that represent
$\H^n(\Cur\sl_2, M_U)$ for any $\sl_2$-module $U$. 
Namely, $\H^n(\Cur\sl_2, M_U) \simeq
\Hom_{\sl_2}(\Sym^n \sl_2 / (h^2-4ef), U)
\oplus \Hom_{\sl_2}(\Sym^{n-3} \sl_2 / (h^2-4ef), U)$.
The cocycle $\al\in C^n(\Cur\sl_2, M_U)$ that corresponds 
to $(\varphi_n,\varphi_{n-3})$ is
\begin{multline*}
\al_{\la_1,\dots,\la_n}(a_1,\dots,a_n)
\\
\begin{split}
&= \Pi(\la_1,\dots,\la_n) \,
\varphi_n(a_1,\dots, a_n)
\\
&+ \sum_{1\le i<j<k\le n} c_3(a_i,a_j,a_k) \,
\Pi(\la_1,\dots,\what\la_i,\dots,\what\la_j,\dots,\what\la_k,\dots,\la_n)
\end{split}
\\
\times
\varphi_{n-3}
(a_1,\dots,\what a_i,\dots,\what a_j,\dots,\what a_k,\dots, a_n),
\end{multline*}
where
$\Pi(\la_1,\dots,\la_n) = \la_1\dotsm\la_n \prod_{1\le r<s\le n}(\la_r-\la_s)$
and $c_3(a_1,a_2,a_3) = (a_1\wedge a_2\wedge a_3)/(e\wedge f\wedge h)$
is the generator of $\H^3(\sl_2) \simeq \CC$.
\end{ex}

\begin{remark}
Corollary~\ref{h1,2cur} in light of Theorem~\ref{exts} implies the
following explicit description of the two-cocycles $c_\la(a,b)$
corresponding to abelian extensions 
\[
0 \to M_U \to \wti A \to A \to 0
\]
of a current conformal algebra $A = \Cur \gtg$ by a
current module $M_U$. (See the proof of  Theorem~\ref{exts}, Part 4
for the notation.)

When $\gtg \ne \sl_2$, abelian extensions are parameterized by
elements $\varphi\in \Hom_\g(\bigwedge^2\g / \g, U)$ and the
corresponding cocycle is $c_\la(a,b) \linebreak[1] = \linebreak[0]
\la(\d \nolinebreak + \nolinebreak \la)\varphi(a,b)$.

When $\gtg = \sl_2$, abelian extensions are parameterized by 
elements 
$\varphi\in \Hom_{\sl_2}(\Sym^2 \sl_2 / (h^2-4ef), U)
= \Hom_{\sl_2}(V(4), U)$ 
and
$c_\la(a,b) = \la(\d+\la) (\d+2\la)\varphi(a,b)$.
\end{remark}

%%%%%%%%%%%%%%%%%%%%%%%%
\section{Hochschild, cyclic, and Leibniz \coh}\label{shc}
%%%%%%%%%%%%%%%%%%%%%%%%
\subsection{Hochschild \coh}
\label{hochschild}

We can similarly define the notion of Hochschild \coh\ by considering
the following analogues of the basic and reduced complexes for an
associative conformal algebra $A$ and a conformal bimodule $M$ over
it, see Definition~\ref{df:1.4}.
\begin{df}
%\label{cochains}
A \emph{Hochschild $n$-cochain $(n \in \nz_+)$ of an associative
conformal algebra $A$ with coefficients in a conformal bimodule $M$
over it} is a $\nc$-linear operator
\begin{align*}
\gamma: A^{\otimes n} & \to M[\lambda_1, \dots, \lambda_n]\\
a_1 \otimes \dots \otimes a_n & \mapsto \gamma_{\lambda_1, \dots,
\lambda_n}(a_1, \dots, a_n),
\end{align*}
satisfying the following condition:
\begin{description}
\item[Conformal antilinearity] $\gamma_{\lambda_1, \dots, \lambda_n}
(a_1, \dots, \partial a_i, \dots, a_n)$ \\  $= - \lambda_i
\gamma_{\lambda_1, \dots, \lambda_n} (a_1, \dots, a_i, \dots, a_n)$
for all $i$.
%\item[Locality] $\gamma_{\lambda_1, \dots, \lambda_n}(a_1, \dots, a_n)$
%depends polynomially on the $\lambda_i$'s for any fixed $a_1, \dots,
%a_n \in A$.
\end{description}
\end{df}
The differential $d$ of a cochain $\gamma$ is defined as follows:
\begin{align*}
&(d \gamma)_{\lambda_1, \dots, \lambda_{n+1}} (a_1,\dots,a_{n+1})
\\
& = {a_1}_{\la_1} \gamma_{\lambda_2, \dots, \lambda_{n+1}}
(a_2,\dots,a_{n+1})
\\
& \begin{split}
+ \sum^{n}_{i = 1} (-1)^{i} \gamma_{\lambda_1, \dots, \lambda_{i-1},
\lambda_i+\lambda_{i+1}, \lambda_{i+2}, \dots, \lambda_{n+1}} (a_1, \dots,
a_{i-1}, \, & {a_i}_{\lambda_i} a_{i+1}, 
\\
& a_{i+2}, \dots, a_{n+1})
\end{split}
\\
& \phantom{+}
+ (-1)^{n+1} \gamma_{\lambda_1, \dots, \lambda_{n}}
(a_1,\dots,a_{n})_{-\d - \la_{n+1}} a_{n+1}.
\end{align*}
One can verify that the operator $d$ preserves the space of cochains
and $d^2 = 0$.  The cochains of an associative conformal algebra $A$
with coefficients in a bimodule $M$ form a complex $\ct^\bullet =
\ct^\bullet (A,M)$, called the \emph{basic Hochschild complex}. As in
the Lie conformal algebra case, $\ct^\bullet (A,M)$ carries the
structure of a (left) $\nc[\partial]$-module:
\begin{equation}
\label{d-module}
(\partial\cdot \gamma)_{\lambda_1, \dots, \lambda_n}(a_1, \dots, a_n) 
= \Bigl( \partial_M + \sum_{i=1}^n \lambda_i \Bigr) 
\gamma_{\lambda_1, \dots, \lambda_n}(a_1,\dots, a_n).
\end{equation}
A straightforward computation shows that $d$ commutes with $\d$. The
quotient complex
\[
C^\bullet (A,M) = \ct^\bullet (A,M)/\partial \ct^\bullet (A,M)
\]
is called the \emph{reduced Hochschild complex}, and its \coh\ is
called the \emph{reduced Hochschild \coh} $\H^\bullet(A,M)$, as
opposed to the \emph{basic Hochschild \coh} $\wti \H^\bullet(A,M)$,
which is the
\coh\ of the basic Hochschild complex $\ct^\bullet$. Low-degree
Hochschild \coh\ groups can be interpreted along the lines of
Section~\ref{s4}, \emph{e.g}., 
$\wti \H^0 (A,M) = \{ m \in M \; | \;
a_\la m = m_{-\d-\la} a \;\;\forall a\in A\}$. 
%The Hochschild
%\coh\ is harder to compute, as usually --- that was the main reason of
%our resorting to the conformal analogue of Chevalley--Eilenberg \coh\
%in this paper.

%%%%%%%%%%%%%%%%%%%%%%%%%%%
\begin{remark}\label{ran6.1,6.2}
One has obvious analogues of Theorems \ref{ec1}, \ref{ec2}, and
Proposition~\ref{hcurgmu} for Hochschild cohomology. 
\end{remark}
%%%%%%%%%%%%%%%%%%%%%%%%%%%

For a current conformal algebra $\Cur A$, where $A$ is 
a $\nc$-algebra, the reduced Hochschild \coh\ 
$\H^\bullet (\Cur A, \Cur A) \simeq \H^\bullet(A[t],A[t])$,
by the analogue of Proposition~\ref{hcurgmu}.
By the Hochschild--Kostant--Rosenberg Theorem \cite{HKR},
when $A$ is the algebra of regular
functions on an affine nonsingular scheme $\Spec A$ over $\nc$, 
the latter \coh\ is isomorphic to the space of polyvector fields
$\bigwedge_A T_A \otimes_\nc \bigwedge_{\nc[t]} \nc[t]
\d_t$ on the product $\Spec A \times \na^1$,
where $T_A = \Der (A,A)$ is the left module of vector fields on
$\Spec A$.

%%%%%%%%%%%%%%%%%%%%%%%%%%%
\begin{remark}
\begin{sloppypar}
For a commutative associative conformal algebra \cite{K2}, one
can define the analogue of the Harrison \coh\ by placing the symmetry
condition on Hochschild cochains. This \coh\ is the closest analogue
of the one introduced in \cite{KV} in the context of vertex algebras.
\end{sloppypar}
\end{remark}
%%%%%%%%%%%%%%%%%%%%%%%%%%%
\subsection{Cyclic \coh}\label{scyclcoh}

In this section we define an analogue of cyclic \coh, see
\cite{connes,loday,tsygan}, for an associative conformal algebra $A$. Define
its \emph{basic cyclic \coh} $\wti{\HC}^\bullet (A)$ as the \coh\ of
the complex $\wti\Cc^\bullet$, where $\wti\Cc^n$, $n\in\ZZ_+$, is the
space of $\nc$-linear operators
\begin{align*}
\gamma: A^{\otimes (n+1)} & \to \nc[\lambda_0, \dots, \lambda_{n}]\\
a_0 \otimes \dots \otimes a_n & \mapsto \gamma_{\lambda_0, \dots,
\lambda_n}(a_0, \dots, a_n),
\end{align*}
satisfying the following conditions:
\begin{description}
\item[Conformal antilinearity] $\gamma_{\lambda_0, \dots, \lambda_n}
(a_0, \dots, \partial a_i, \dots, a_n)$ \\  $= - \lambda_i
\gamma_{\lambda_0, \dots, \lambda_n} (a_0, \dots, a_i, \dots, a_n)$;
%\item[Locality] $\gamma_{\lambda_0, \dots, \lambda_n}(a_0, \dots, a_n)$
%depends polynomially on the $\lambda_i$'s for any fixed $a_0, \dots,
%a_n \in A$;
\item[Cyclic invariance] $\gamma_{\la_1, \dots, \la_n, \la_0}
(a_1, \dots, a_n, a_0)$\\ $= \linebreak[1] (-1)^n
\gamma_{\la_0, \dots, \la_n} (a_0,
\linebreak[0] \dots, a_n)$.
\end{description}
The differential $d$ of a cochain $\gamma$ is defined as follows:
\begin{align*}
&(d \gamma)_{\lambda_0, \dots, \lambda_{n+1}} (a_0,\dots,a_{n+1})
\\
& \begin{split}
 = \sum^{n}_{i = 0} (-1)^{i} \gamma_{\lambda_0, \dots, \lambda_{i-1},
\lambda_i+\lambda_{i+1}, \lambda_{i+2}, \dots, \lambda_{n+1}} (a_0, \dots,
a_{i-1}, \, & {a_i}_{\lambda_i} a_{i+1}, 
\\
& a_{i+2}, \dots, a_{n+1})
\end{split}
\\
& \phantom{+}
+ (-1)^{n+1} \gamma_{\lambda_{n+1}+\la_0, \dots, \lambda_{n}}
({a_{n+1}}_{\la_{n+1}} a_0,\dots,a_{n}).
\end{align*}

The \emph{reduced cyclic \coh} $\HC^\bullet (A)$ may be defined as the
\coh\ of the quotient complex by the action of $\d$, as in the
Hochschild case. 

\subsection{Leibniz \coh}
\label{leibniz}

Nonlocal collections of formal distributions lead to the notion of a
Leibniz conformal algebra, see Section~\ref{s1}:
\begin{df}
A \emph{Leibniz conformal algebra} is a $\CC [\partial ]$-module $A$
endowed with a $\lambda$-bracket $[a_{\lambda} b]$ which defines a
conformally sesquilinear map $A \otimes A \to A [[ \lambda ]]$
satisfying the Jacobi identity as in Definition~\ref{df:1.1}.
\end{df}
The difference from Definition~\ref{df:1.1} of a Lie conformal algebra
is that the skew-symmetry axiom is omitted and formal power series in
$\la$ are allowed. For a Leibniz conformal algebra $A$, the definition
of a (\emph{left$)$ module $M$ over it} is the same as that for Lie
conformal algebras, see Definition~\ref{df:1.2}. The \emph{space
$C^n(A,M)$ of $n$-cochains of a Leibniz algebra $A$ with values in a
module} $M$ is the space of $\nc$-linear operators
\begin{align*}
\gamma: A^{\otimes n} & \to M[[\lambda_1, \dots, \lambda_n]]\\
a_1 \otimes \dots \otimes a_n & \mapsto \gamma_{\lambda_1, \dots,
\lambda_n}(a_1, \dots, a_n),
\end{align*}
which are conformally antilinear:
\begin{equation*}
\gamma_{\lambda_1, \dots, \lambda_n}
(a_1, \dots, \partial a_i, \dots, a_n)
= - \lambda_i
\gamma_{\lambda_1, \dots, \lambda_n} (a_1, \dots, a_i, \dots, a_n).
\end{equation*}
The differential $d$ of a cochain $\gamma$ is defined as follows:
\begin{multline*}
(d\gamma)_{\lambda_1, \dots, \lambda_{n+1}} (a_1,\dots,a_{n+1})
\\
\begin{split}
& = \sum_{i=1}^{n+1} (-1)^{i+1}
{a_i}_{\lambda_i} \gamma_{\lambda_1, \dots, \widehat\lambda_i,
\dots, \lambda_{n+1}} (a_1,\dots, \widehat a_i,\dots, a_{n+1}) 
\\
&
\begin{split}
 \phantom{=} + \sum_{1 \le i < j \le n+1} (-1)^{i} \gamma_{\lambda_1, \dots,
\widehat\lambda_i, \dots,
\lambda_{j-1}, \lambda_i+\lambda_j, \la_{j+1}, \dots,
\lambda_{n+1}} (a_1, \dots, \widehat a_i,\\
\dots, a_{j-1}, [{a_i} _{\lambda_i} a_j], a_{j+1}, \dots, a_{n+1}),
\end{split}
\end{split}
\end{multline*}
where $\gamma$ is extended linearly over the polynomials in
$\lambda_i$.  One can verify that the operator $d$ preserves the space
of cochains and $d^2 = 0$.  The $n$-cochains, $n \in \nz_+$, of a
Leibniz conformal algebra $A$ with coefficients in a module $M$ form
a complex $\ct^\bullet = \ct^\bullet (A,M)$, called the \emph{basic
Leibniz complex}. 

\begin{sloppypar}
Equation~\eqref{d-module} defines the structure of a left
$\nc[\partial]$-module on $\ct^\bullet (A,M)$, which commutes with
$d$. The quotient complex
\[
C^\bullet (A,M) = \ct^\bullet (A,M)/\partial \ct^\bullet (A,M)
\]
is called the \emph{reduced Leibniz complex}. Its \coh\ is called the
\emph{reduced Leibniz \coh} $\H^\bullet(A,M)$, as opposed to the
\emph{basic Leibniz \coh} $\wti \H^\bullet(A,M)$, which is the \coh\
of the basic Leibniz complex $\ct^\bullet$. These are conformal
analogues of cohomology of Leibniz algebras, see
\cite{cuvier,loday,loday-leib}.
\end{sloppypar}

%%%%%%%%%%%%%%%%%%%%%%%%%%%%%%%%%%%%%%%%%%%%%%%%%%%%%%
\section{Generalization to conformal algebras in 
higher dimensions}
\label{s12}
%%%%%%%%%%%%%%%%%%%%%%%%%%%%%%%%%%%%%%%%%%%%%%%%%%%%%%
The theory of conformal algebras, their representations and cohomology
has a straightforward generalization to the case when $\la$ is a
vector.

Let us fix a natural number $r$.  We replace a single indeterminate
$\la$ by the vector $\vec\la = (\la_1,\dots,\la_r)$ and $\d$ by
$\vec\d = (\d_1,\dots,\d_r)$, and use the multi-index notation like
$\vec\la^{(\vec m)} = \la_1^{(m_1)}\dotsm\la_r^{(m_r)}$ for $\vec
m\in\Zset^r$, $\de (\vec z- \vec w) = \prod_i \de (z_i-w_i)$, etc.
Then everything from Sections 1--6 and 9 holds.

Examples of conformal algebras in $r$ indeterminates are provided by
$r$-dimensional current algebras, cf.\ Examples \ref{ex:1.1} and \ref{dncurr}.
Other important examples are the Cartan algebras of vector fields.
The structure theory of higher dimensional conformal algebras,
including a classification of the simple ones, is currently being developed
\cite{BK}.

%%%%%%%%%%%%%%%%%%%%%%%%%%%
\begin{ex}\label{exwr}
The Lie algebra 
$W_r = \Der\Cset[x_1,x_1^{-1},\dots,x_r,x_r^{-1}]$
is spanned by the coefficients of the formal distributions
\begin{equation*}
L^i(\vec z) = -\de(\vec z - \vec x)\d_{x_i} 
\; \Bigl(=-\sum_{\vec m\in\Zset^r} x_1^{m_1}\dotsm x_r^{m_r}  \d_{x_i}
z_1^{-m_1-1}\dotsm z_r^{-m_r-1} \Bigr).
\end{equation*}
They are pairwise local, since
\begin{multline*}
[L^i(\vec z),L^j(\vec w)] 
= \d_{w_i}\bigl(L^j(\vec w) \de(\vec z-\vec w)\bigr)
- \d_{z_j}\bigl(L^i(\vec w) \de(\vec z-\vec w)\bigr)
\\
= \d_{w_i}L^j(\vec w) \de(\vec z-\vec w)
+ L^j(\vec w) \d_{w_i}\de(\vec z-\vec w)
+ L^i(\vec w) \d_{w_j}\de(\vec z-\vec w).
\end{multline*}
The corresponding conformal algebra is 
$A=\bigoplus_{i=1}^r \Cset[\vec\d] L^i$
with $\vec\la$-brackets
\begin{equation}\label{laprodwr}
[{L^i}_{\vec\la} L^j] = \d_i L^j + \la_i L^j + \la_j L^i.
\end{equation}
Its annihilation algebra is 
${W_r}_- = \Der\Cset[x_1,\dots,x_r]$.
For $r=1$ $A$ is the Virasoro conformal algebra $\vir$, 
see Example~\ref{ex:1.2}.

By Corollary \ref{cec1}, the cohomology of ${W_r}_-$ with trivial
coefficients is the same as the cohomology of the complex 
$\wti C^\bullet(A,\Cset)$. The latter can be described as follows.
Let $V$ be the vector space $\bigoplus_{i=1}^r \Cset L^i$.
Every cochain $\al\in \wti C^n(A,\Cset)$ is uniquely determined by its values
on $V^{\tt n}$:
\begin{equation*}
\al\colon V^{\tt n}\to\Cset[\vec\la_1,\dots,\vec\la_n],
\quad
L^{k_1}\tt\dots\tt L^{k_n} \mapsto 
\al_{\vec\la_1,\dots,\vec\la_n}^{k_1,\dots,k_n}.
\end{equation*}
The differential is given by the formula
\begin{align*}%\label{difwr}
(d\al)_{\vec\la_1,\dots,\vec\la_{n+1}}^{k_1,\dots,k_{n+1}}
&= \sum^{n+1}_{\substack{i,j = 1 \\ i < j}} (-1)^{i+j}
\la_{i,k_j} \al_{\vec\la_i+\vec\la_j,\vec\la_1,\dots,
\widehat{\vec\la}_i,\dots,\widehat{\vec\la}_j,\dots,\vec\la_{n+1}}^{\;\;\;\;\;
\;\;
k_i,k_1,\dots,\widehat{k_i},\dots,\widehat{k_j},\dots,k_{n+1}}
\\
&- \sum^{n+1}_{\substack{i,j = 1 \\ i < j}} (-1)^{i+j}
\la_{j,k_i} \al_{\vec\la_i+\vec\la_j,\vec\la_1,\dots,
\widehat{\vec\la}_i,\dots,\widehat{\vec\la}_j,\dots,\vec\la_{n+1}}^{\;\;\;\;\;
\;\;
k_j,k_1,\dots,\widehat{k_i},\dots,\widehat{k_j},\dots,k_{n+1}} ,
\end{align*}
where $\la_{i,k}$ is the $k$th coordinate of the vector ${\vec\la}_i$.

The cohomology of the Lie algebra ${W_r}_-$ with trivial
coefficients was computed by Gelfand and Fuchs \cite{GF2} 
(see also \cite[$\S 2.2.2$]{F}).
\end{ex}
%%%%%%%%%%%%%%%%%%%%%%%%%%%
\begin{ex}\label{exsr}
The subalgebra of divergence $0$ derivations is a formal distribution
subalgebra of $W_r$. The corresponding conformal algebra is the
following subalgebra of the algebra in Example~\ref{exwr}:
%\begin{equation*}
$\{ \sum_i P_i(\vec\d)L^i \,|\linebreak[1] \sum_i P_i(\vec\d) \d_i = 0 \}$.
%\end{equation*}
\end{ex}
%%%%%%%%%%%%%%%%%%%%%%%%%%%
\begin{ex}\label{exhr}
The subalgebra $H_r$, $r=2s$, of Hamiltonian derivations is a formal 
distribution subalgebra of $W_r$. The corresponding conformal algebra is
of rank one: $A=\Cset[\vec\d]L$ with $\vec\la$-bracket
\begin{equation*}
[L_{\vec\la}L] = \sum_{i=1}^s (\la_{s+i} \d_i L - \la_i \d_{s+i} L).
\end{equation*}
Its annihilation algebra ${H_r}_-$ is the Lie algebra of Hamiltonian 
derivations of $\Cset[x_1,\dots,x_r]$.

\begin{sloppypar}
The $n$th term of the complex $\wti C^\bullet(A,\Cset)$, whose cohomology is 
$\H^\bullet({H_r}_-)$,
can be identified with the space of skew-symmetric polynomials in
$\vec\la_1,\dots,\vec\la_n$. 
The differential is given by the formula
\begin{multline*}%\label{difhr}
(d P)(\vec\la_1,\dots,\vec\la_{n+1})
\\*
= \sum^{n+1}_{\substack{i,j = 1\\i < j}} (-1)^{i+j}
(\vec\la_i | \vec\la_j) 
P(\vec\la_i+\vec\la_j,\vec\la_1,\dots,
\widehat{\vec\la}_i,\dots,\widehat{\vec\la}_j,\dots,\vec\la_{n+1})
\end{multline*}
where 
$(\vec\la | \vec\mu) = \sum_{k=1}^s (\la_k \mu_{s+k} - \la_{s+k} \mu_k)$.
\end{sloppypar}

For $r=2$ this complex has been known for quite a long time, 
but the computation
of its cohomology is still an open problem (see \cite[$\S 2.2.7$]{F}).
\end{ex}
%%%%%%%%%%%%%%%%%%%%%%%%%%%
\begin{ex}\label{exkr}
The subalgebra $K_r$, $r=2s+1$, of contact derivations is also a formal 
distribution subalgebra of $W_r$, but the corresponding conformal algebra is
of infinite rank. It is better viewed as a Lie${}^*$ algebra of rank $1$,
see Section~\ref{s13}.
\end{ex}

%%%%%%%%%%%%%%%%%%%%%%%%%%%%%%%%%%%%%%%%%%%%%%%%%%%%%%
\section{Higher differentials}\label{s11}
%%%%%%%%%%%%%%%%%%%%%%%%%%%%%%%%%%%%%%%%%%%%%%%%%%%%%%
For the computation of the cohomology with non-trivial coefficients
of the Lie algebras of vector fields,
it is useful to know the cohomology of their subalgebras
of vector fields which have a zero of certain order at the origin
(see \cite{F}). The argument of Theorem \ref{ec1} can be generalized
to give a complex which produces this cohomology.

Let $A$ be a conformal algebra in $r$ indeterminates which is a free
$\Cset[\vec\d]$-module: $A=\bigoplus_{i\in I} \Cset[\vec\d] L^i$.  For
fixed $\vec N\in\Zset_+^r$, we define $\g_{\vec N}\equiv (\Lie
A)_{\vec N}$ to be the subspace of the annihilation algebra
$\g_-=(\Lie A)_-$, spanned by $L^i_{\vec m}$, $i\in I$, $\vec m\ge
\vec N$ (meaning that $m_i \ge N_i$ for each $i$).  We are interested
in the case when $\g_{\vec N}$ is a {\em Lie subalgebra\/} of $\g_-$.
Note that this is always true when the entries of $\vec N$ are large
enough.  Indeed, we can write
\begin{equation}\label{strconst}
[L^i_{\vec\la}L^j] = \sum_{k\in I} C_{ij}^k(\vec\la,\vec\d) L^k
\end{equation}
for some uniquely determined polynomials $C_{ij}^k$. Then
\begin{equation}\label{llalmu}
[L^i_{\vec\la}, L^j_{\vec\mu}] 
= \sum_{k\in I} C_{ij}^k(\vec\la, -\vec\la-\vec\mu) L^k_{\vec\la+\vec\mu}.
\end{equation}
It follows that for large $\vec N$ the commutator 
$[\d_{\vec\la}^{\vec N} L^i_{\vec\la}, 
  \d_{\vec\mu}^{\vec N} L^j_{\vec\mu}]$
can be expressed in terms of 
$\d_{\vec\la+\vec\mu}^{\vec N} L^k_{\vec\la+\vec\mu}$.
Since
$\d_{\vec\la}^{\vec N} L^i_{\vec\la}
= \sum_{\vec m\ge \vec N} L^i_{\vec m} {\vec\la}^{(\vec m - \vec N)}$,
this shows that $\g_{\vec N}$ is a Lie subalgebra of $\g_-$.

Let $M$ be a module over the conformal algebra $A$.
Then $M$ is a $\g_-$-module and hence also a $\g_{\vec N}$-module.
Let $V$ be the vector space $\bigoplus_{i\in I} \Cset L^i$.
As in Section~\ref{s8}, the $n$th term of the complex 
$C^\bullet(\g_{\vec N}, M)$ can be identified with the space
of linear maps
\begin{equation}\label{cN}
\begin{split}
\al\colon V^{\tt n} & \to\Cset[\vec\la_1,\dots,\vec\la_n] \tt_\Cset M,
\\
L^{k_1}\tt\dotsm\tt L^{k_n} & \mapsto 
\al_{\vec\la_1,\dots,\vec\la_n}^{k_1,\dots,k_n},
\end{split}
\end{equation}
which are skew-symmetric with respect to simultaneous permutations
of $k_i$'s and $\vec\la_i$'s. Using \eqref{llalmu}, one can easily
write its differential $d_{\vec N}$.
%%%%%%%%%%%%%%%%%%
\begin{ex}\label{dnwr}
Let $A$ be the conformal algebra
associated to the Lie algebra $W_r$ of vector fields (Example \ref{exwr}).
Then for $\vec N\in\Zset_+^r$, $W_{r,\vec N}$ is the Lie algebra of
vector fields 
$\sum P_i(\vec x)\d_{x_i}$
such that all $P_i(\vec x)$ are divisible by ${\vec x}^{\vec N}$.
Equation \eqref{laprodwr} implies
\begin{equation*}
[\d_{\vec\la}^{\vec N} L^i_{\vec\la}, 
  \d_{\vec\mu}^{\vec N} L^j_{\vec\mu}]
= \d_{\vec\la}^{\vec N} \la_j 
\d_{\vec\la+\vec\mu}^{\vec N} L^i_{\vec\la+\vec\mu}
- \d_{\vec\mu}^{\vec N} \mu_i 
\d_{\vec\la+\vec\mu}^{\vec N} L^j_{\vec\la+\vec\mu}.
\end{equation*}
The differential $d_{\vec N}$ of the complex \eqref{cN} is given by the
formula
\begin{align*}
(d_{\vec N} \al)_{\vec\la_1,\dots,\vec\la_{n+1}}^{k_1,\dots,k_{n+1}}
&= \sum_{i=1}^{n+1} (-1)^{i+1}\,\d_{\vec\la_i}^{\vec N} {L^{k_i}}_{\vec\la_i}
\al_{\vec\la_1,\dots,\widehat{\vec\la}_i,\dots,
\vec\la_{n+1}}^{k_1,\dots,\widehat{k_i},\dots,k_{n+1}}
\\
&+ \sum^{n+1}_{\substack{i,j = 1 \\ i < j}} (-1)^{i+j}\,
\d_{\vec\la_i}^{\vec N}\la_{i,k_j} \al_{\vec\la_i+\vec\la_j,\vec\la_1,\dots,
\widehat{\vec\la}_i,\dots,\widehat{\vec\la}_j,\dots,\vec\la_{n+1}}^{\;\;\;\;\;
\;\;
k_i,k_1,\dots,\widehat{k_i},\dots,\widehat{k_j},\dots,k_{n+1}}
\\
&- \sum^{n+1}_{\substack{i,j = 1 \\ i < j}} (-1)^{i+j}\,
\d_{\vec\la_j}^{\vec N}\la_{j,k_i} \al_{\vec\la_i+\vec\la_j,\vec\la_1,\dots,
\widehat{\vec\la}_i,\dots,\widehat{\vec\la}_j,\dots,\vec\la_{n+1}}^{\;\;\;\;\;
\;\;
k_j,k_1,\dots,\widehat{k_i},\dots,\widehat{k_j},\dots,k_{n+1}}
.
\end{align*}
\end{ex}
%%%%%%%%%%%%%%%%%%
\begin{ex}\label{dncurr}
Let $\g$ be a Lie algebra. Then the current algebra
$\tilde\g = \g\tt_\Cset \Cset[x_1,x_1^{-1},\dots,x_r,x_r^{-1}]$
is spanned by the coefficients of the pairwise local formal distributions
$a(\vec z) := a\tt\de(\vec x - \vec z)$, $a\in\g$. They satisfy
$[a(\vec z), b(\vec w)] = [a,b](\vec w)\de(\vec z - \vec w)$.
The corresponding conformal algebra is $A=\Cset[\vec\d]\tt_\Cset\g$
with $\vec\la$-brackets determined by
\begin{equation*}
[a_{\vec\la}b] = [a,b] \qquad\text{ for }\;\; a,b\in\g.
\end{equation*}
The annihilation algebra of $A$ is 
$\tilde\g_- = \g\tt_\Cset \Cset[\vec x]$
and for $\vec N\in\Zset_+^r$ 
$\tilde\g_{\vec N} = \g\tt_\Cset \Cset[\vec x] {\vec x}^{\vec N}$.
Now $C^n(\tilde\g_{\vec N}, M)$ consists of all 
\begin{align*}
\al\colon \g^{\tt n}& \to M [\vec\la_1,\dots,\vec\la_n] ,
\\
a_1\tt\dotsm\tt a_n & \mapsto 
\al_{\vec\la_1,\dots,\vec\la_n}(a_1,\dots,a_n),
\end{align*}
skew-symmetric with respect to simultaneous permutations
of $a_i$'s and $\vec\la_i$'s. 
The differential $d_{\vec N}$ is given by
\begin{multline*}
(d_{\vec N} \al)_{\vec\la_1,\dots,\vec\la_{n+1}} (a_1,\dots,a_{n+1})
\\*
\begin{split}
&= \sum_{i=1}^{n+1} (-1)^{i+1}\,\d_{\vec\la_i}^{\vec N} {a_i}_{\vec\la_i}
\al_{\vec\la_1,\dots,\widehat{\vec\la}_i,\dots,
\vec\la_{n+1}} (a_1,\dots,\widehat{a_i},\dots,a_{n+1})
\\
&
\begin{split}
+ \sum^{n+1}_{\substack{i,j = 1 \\ i < j}} (-1)^{i+j}\,
\d_{\vec\la_i}^{\vec N} \al_{\vec\la_i+\vec\la_j,\vec\la_1,\dots,
\widehat{\vec\la}_i,\dots,\widehat{\vec\la}_j,\dots,\vec\la_{n+1}}
([a_i,a_j],a_1,\dots,\widehat{a_i},
\\
\dots,\widehat{a_j},\dots,a_{n+1}).
\end{split}
\end{split}
\end{multline*}
\end{ex}
%%%%%%%%%%%%%%%%
\begin{remark}
It is easy to see that in the above examples the differentials satisfy
$d_{\vec N}d_{\vec N'} + d_{\vec N'}d_{\vec N} = 0$.
\end{remark}
%%%%%%%%%%%%%%%%

%%%%%%%%%%%%%%%%%%%%%%%%%%%%%%%%%%%%%%%%%%%%%%%%%%%%%%
\section{Relation to Lie${}^*$ algebras}\label{s13}
%%%%%%%%%%%%%%%%%%%%%%%%%%%%%%%%%%%%%%%%%%%%%%%%%%%%%%
The theory of conformal algebras is in many ways analogous to
the theory of Lie algebras. The reason is that in fact conformal algebras
can be considered as Lie algebras in certain pseudo-tensor categories,
instead of the category of vector spaces. 
A pseudo-tensor category \cite{BD} is a category equipped with 
``polylinear maps'' and a way to compose them. This is enough
to define the notions of Lie algebra, representations, cohomology,
etc.

As an example, consider first the category $\Vec$
of vector spaces (over $\Cset$).
For a finite non-empty set $I$ and
a collection of vector spaces $\{L_i\}_{i\in I}$, $M$,
we can define {\em polylinear maps\/} 
{}from $\{L_i\}_{i\in I}$ to $M$:
\begin{equation*}%\label{plin}
P_I(\{L_i\}_{i\in I}, M) 
:= \Hom(\otimes_{i\in I} L_i, M).
\end{equation*}
This is a vector space with an action of the symmetric group
$\symm_I$ on it.

For any surjection of finite sets $J\overset{\pi}{\surjto} I$
and a collection $\{K_j\}_{j\in J}$,
we have the obvious compositions of polylinear maps
\begin{align}\label{com1}
&P_I(\{L_i\}_{i\in I}, M) 
\otimes \bigotimes_{i\in I}
P_{J_i}(\{K_j\}_{j\in J_i}, L_i) 
\to P_J(\{K_j\}_{j\in J},M),    \\
\label{com2}
& \phi \times \{\psi_i\}_{i\in I} \mapsto 
\phi\circ(\otimes_{i\in I} \psi_i) \equiv \phi(\{\psi_i\}_{i\in I}),
\end{align}
where $J_i := \pi^{-1}(i)$ for $i\in I$.

The compositions have the following properties:
\begin{description}
\item[Associativity] 
\begin{sloppypar}
If $H\surjto J$, $\{F_h\}_{h\in H}$ is a family of objects and
$\chi_j\in P_{H_j}(\{F_h\}_{h\in H_j}, K_j)$, then
$\phi\bigl(\bigl\{\psi_i(\{\chi_j\}_{j\in J_i})\bigr\}_{i\in I}\bigr)
\linebreak[1] = \linebreak[0]
\bigl(\phi(\{\psi_i\}_{i\in I})\bigr)(\{\chi_j\}_{j\in J}) 
\linebreak[0] \in P_H(\{F_h\}_{h\in H}, M)$.
\end{sloppypar}

\item[Unit] 
For any object $M$ there is an element $\id_M\in P_1(\{M\},M)$ such
that for any $\phi\in P_I(\{L_i\}_{i\in I}, M)$ one has
$\id_M(\phi)=\phi(\{\id_{L_i}\}_{i\in I})=\phi$.

\item[Equivariance] 
The compositions \eqref{com1}
are equivariant with respect to the natural action of the symmetric group.
\end{description}
%%%%%%%%%%%%%%%%%%%%%%%%%%%
\begin{df}\label{dptc} \cite{BD}.
A {\em pseudo-tensor category\/} is a class of objects $\M$ 
together with
vector spaces $P_I(\{L_i\}_{i\in I}, M)$
on which the symmetric group $\symm_I$ acts,
and composition maps \eqref{com1},
satisfying the above three properties.
\end{df}
%%%%%%%%%%%%%%%%%%%%%%%%%%%
\begin{remark}\label{rpt}
For a pseudo-tensor category $\M$ and objects $L,M\in\M$, let
$\Hom(L,M)=P_1(\{L\},M)$. This gives 
a structure of an ordinary (additive) category on $\M$ and all
$P_I$ are functors $(\M^\circ)^I\times\M\to\Vec$.
(Here $\M^\circ$ denotes the dual category of $\M$.)
\end{remark}
%%%%%%%%%%%%%%%%%%%%%%%%%%%
\begin{remark}\label{roper}
The notion of pseudo-tensor category is a straightforward generalization
of the notion of operad. By definition, an {\em operad\/} is a
pseudo-tensor category with only one object.
\end{remark}
%%%%%%%%%%%%%%%%%%%%%%%%%%%

It is instructive to think of 
a polylinear map $\phi\in P_n(\{L_i\}_{i=1}^n, M)$
as an operation with $n$ inputs and $1$ output,
represented by the figure
\begin{equation*}
{\vcenter{\epsfbox{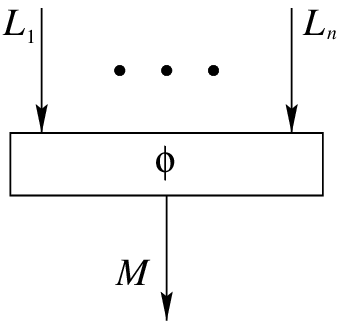}}}
\end{equation*}

%%%%%%%%%%%%%%%%%%%%%%%%%%%
\begin{df}\label{dlie*} 
A {\em Lie algebra in a pseudo-tensor category $\M$\/} 
%(or a {\em Lie${}^*$ algebra})
is an object $A$ and $\mu\in P_2(\{A,A\},A)$ with the following properties.

\begin{description}
\item[Skew-symmetry]
%\begin{equation*}%\label{skew}
$\mu = -\sigma_{12} \, \mu,$
%\end{equation*}
where $\sigma_{12} = (12) \in\symm_2$.

\item[Jacobi identity]
%\begin{equation*}%\label{jac}
$\mu(\mu(\cdot, \cdot), \cdot)=\mu(\cdot, \mu(\cdot, \cdot))-
\sigma_{12} \, \mu(\cdot, \mu(\cdot, \cdot)),$
%\end{equation*}
where now $\sigma_{12} = (12)$ is viewed as an element of $\symm_3$.
\end{description}
\end{df}
%%%%%%%%%%%%%%%%%%%%%%%%%%%

Pictorially, the skew-symmetry and the Jacobi identity for a Lie algebra
$(A,\mu)$ look as follows:
\begin{equation*}
{\vcenter{\epsfbox{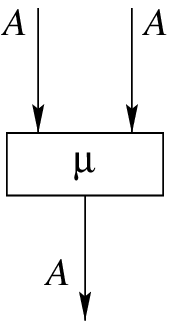}}}
\hspace{-10cm}
= -
{\vcenter{\epsfbox{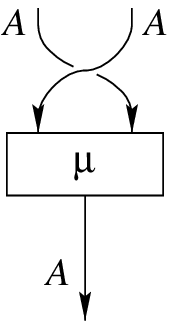}}}
\end{equation*}
\begin{equation*}
{\vcenter{\epsfbox{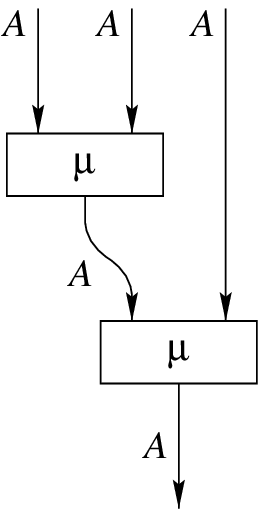}}}
\hspace{-9cm}
=
{\vcenter{\epsfbox{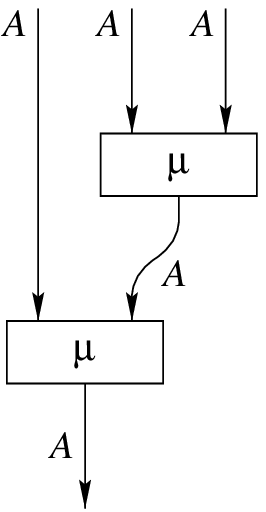}}}
\hspace{-9cm}
+
{\vcenter{\epsfbox{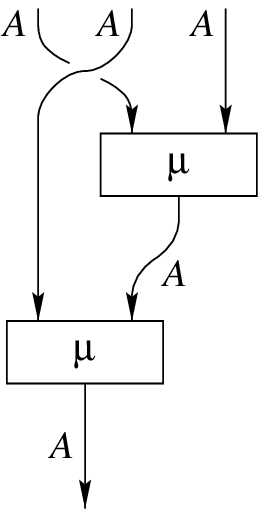}}}
\end{equation*}
%%%%%%%%%%%%%%%%%%%%%%%%%%%
\begin{df}\label{drep}
A {\em representation} of a  Lie algebra $(A,\mu)$ is an object $M$
together with $\rho\in P_2(\{A,M\},M)$ satisfying
\begin{equation*}%\label{rep*}
\rho(\mu(\cdot, \cdot), \cdot)=\rho(\cdot, \rho(\cdot, \cdot))-
\sigma_{12} \, \rho(\cdot, \rho(\cdot, \cdot)).
\end{equation*}
\end{df}
%%%%%%%%%%%%%%%%%%%%%%%%%%%
\begin{df}\label{dcl*}
%\begin{sloppypar}
An {\em $n$-cochain\/} of a Lie algebra $(A,\mu)$ with coefficients in
a module $(M,\rho)$ over it is a polylinear operation $\al
\linebreak[1] \in \linebreak[0] P_n(\{A,\dots,\linebreak[0] A\},M)$
which is skew-symmetric, \emph{i.e.}, satisfying
\begin{equation*}
{\vcenter{\epsfbox{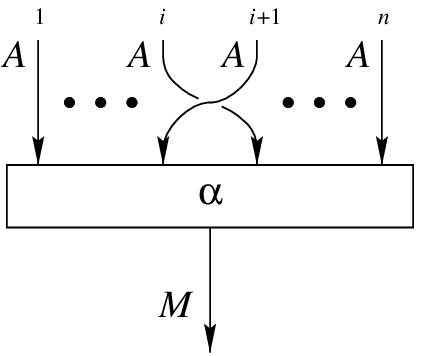}}}
\hspace{-7cm}
= -
{\vcenter{\epsfbox{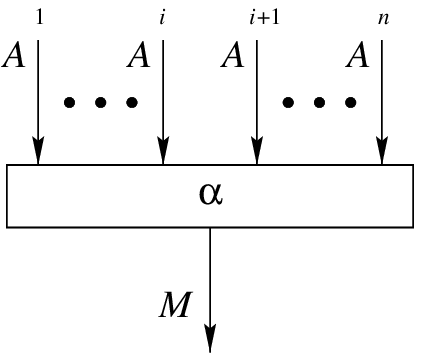}}}
\end{equation*}
for all $i=1,\dots,n$.
%\end{sloppypar}

The differential of a cochain is defined as follows:
{\allowdisplaybreaks{
\begin{align*}
&{\vcenter{\epsfbox{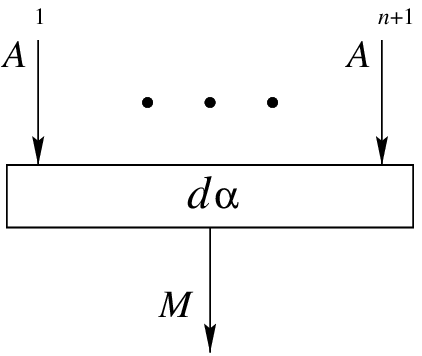}}}
\\
&\hspace{1cm}
= \sum_{1\le i\le n+1} (-1)^{i+1}
{\vcenter{\epsfbox{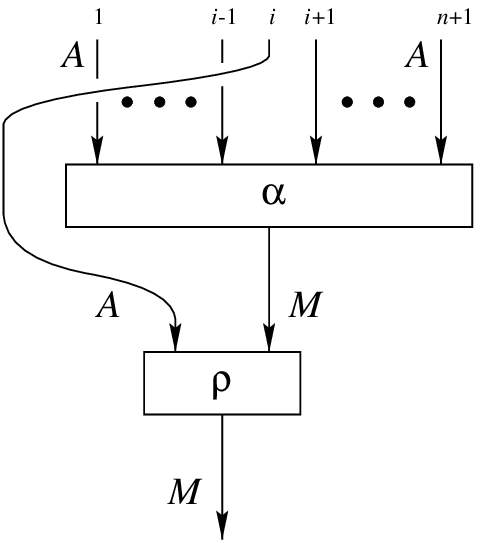}}}
\\
&\hspace{1cm}
+ \sum_{1\le i<j \le n+1} (-1)^{i+j}
{\vcenter{\epsfbox{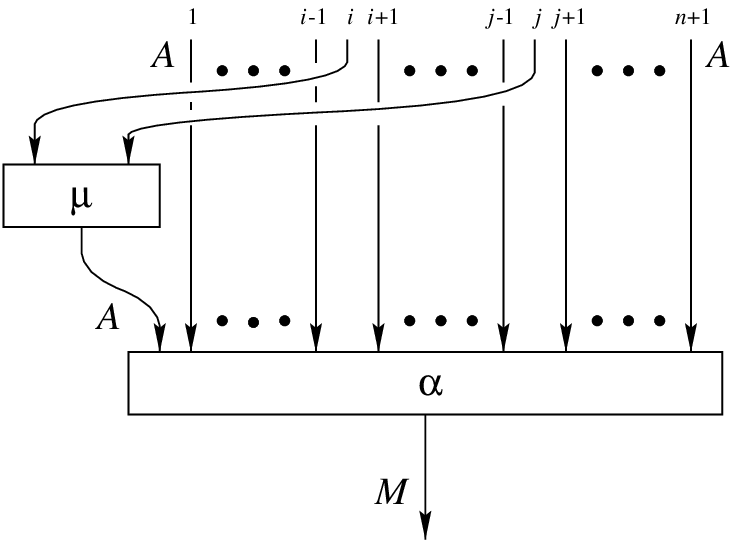}}}
\end{align*}
}}

The same computation as in the ordinary Lie algebra case shows that
$d^2=0$.  The cohomology of the resulting complex is called the 
({\em{reduced}}) 
{\em cohomology of $A$ with coefficients in $M$\/} and is denoted
by $\H^\bullet(A,M)$.
\end{df}
%%%%%%%%%%%%%%%%%%%%%%%%%%%
\begin{remark}\label{rasco*}
One can also define the notions of associative algebra or 
commutative algebra in a pseudo-tensor category, their representations and 
analogues of the Hochschild, cyclic,  or Harrison cohomology.
\end{remark}
%%%%%%%%%%%%%%%%%%%%%%%%%%%
\begin{ex}\label{exuslie}
A Lie algebra in the category of vector spaces $\Vec$ is just
an ordinary Lie algebra. The same is true for representations
and cohomology.
\end{ex}
%%%%%%%%%%%%%%%%%%%%%%%%%%%
\begin{ex}\label{exmld}
Let $D$ be a cocommutative bialgebra with comultiplication $\De$
and counit $\ep$. Then the category $\M^l(D)$ of left $D$-modules
is a symmetric tensor category. Hence, $\M^l(D)$ is a 
pseudo-tensor category with polylinear maps
\begin{equation}\label{mld}
P_I(\{L_i\}_{i\in I}, M) 
:= \Hom_D(\otimes_{i\in I} L_i, M).
\end{equation}
A Lie algebra in the category $\M^l(D)$ is an ordinary Lie algebra which is
also a left $D$-module and such that its bracket is a homomorphism of
$D$-modules.
\end{ex}
%%%%%%%%%%%%%%%%%%%%%%%%%%%
\begin{ex}\label{exm*d}
Let $D$ be as in  Example \ref{exmld}.
We introduce a pseudo-tensor category $\M^*(D)$ with the same objects as 
$\M^l(D)$ but with another pseudo-tensor structure \cite{BD}
\begin{equation}\label{m*d}
P_I(\{L_i\}_{i\in I}, M) 
:= \Hom_{D^{\tt I}} (\boxtimes_{i\in I} L_i, D^{\tt I}\tt_D M).
\end{equation}
Here $\boxtimes_{i\in I}$ is the tensor product functor
$\M^l(D)^I \to \M^l(D^{\tt I})$.
For $J\overset{\pi}{\surjto} I$ the composition of polylinear maps
is defined as follows:
\begin{equation}\label{*com}
\phi\bigl(\{\psi_i\}_{i\in I}\bigr) := \De^{(\pi)}\bigl(\phi\bigr) \,
\circ\bigl(\boxtimes_{i\in I}\psi_i\bigr).
\end{equation}
Here $\De^{(\pi)}$ is the functor 
$\M^l(D^{\tt I})\to\M^l(D^{\tt J})$,
$M\mapsto D^{\tt J} \tt_{D^{\tt I}} M$
where $D^{\tt I}$ acts on $D^{\tt J}$ via the iterated
comultiplication determined by $\pi$.
The symmetric group $\symm_I$ acts on $P_I(\{L_i\}_{i\in I}, M)$
by simultaneously permuting the factors in $\boxtimes_{i\in I} L_i$
and $D^{\tt I}$.
\end{ex}
%%%%%%%%%%%%%%%%%%%%%%%%%%%
\begin{df}\label{dflie*}
A {\em{Lie${}^*$ algebra\/}} is a Lie algebra in the pseudo-tensor category
$\M^*(D)$ defined above.
\end{df}
%%%%%%%%%%%%%%%%%%%%%%%%%%%
The following examples of Lie${}^*$ algebras are important:

1. When $D=\Cset$ we recover Example \ref{exuslie}.

2. For $D=\Cset[\d]$ (with $\De(\d)=\d\tt 1 + 1\tt\d$, $\ep(\d)=0$) we
get exactly the notions of conformal algebras, conformal modules over
them and the reduced cohomology theory introduced in this paper.

3. For $D=\Cset[\d_1,\dots,\d_r]$ we get conformal algebras
in $r$ indeterminates, see Section~\ref{s12}.

4. When $D=\Cset[\Ga]$ is the group algebra of a 
group $\Ga$, one obtains the $\Ga$-conformal algebras
studied in \cite{GK}.

\begin{sloppypar}
5. Let $\Ga$ be a subgroup of $\Cset^*$ and let $D =
\Cset[\d]\rtimes\Cset[\Ga] \linebreak[1] = \linebreak[0]
\bigoplus_{m\in\Zset_+, \al\in\Ga}\,\Cset\, \d^m T_\al$ with
multiplication $T_\al T_\be = T_{\al\be}$, $T_1=1$, $T_\al \d
T_\al^{-1} = \al\d$ and comultiplication $\De(\d)=\d\tt 1 + 1\tt\d$,
$\De(T_\al)=T_\al \tt T_\al$.  Then we get the $\Ga$-conformal
algebras studied in \cite{BK} (cf.\ \cite{K2}).
\end{sloppypar}

6. Let now $D = \Cset[\d]\times F(\Ga)$, where $F(\Ga)$ is the
function algebra of a commutative group $\Ga$.  In other words, $D=
\bigoplus_{m\in\Zset_+, \al\in\Ga}\,\Cset\, \d^m \pi_\al$ with
multiplication $\pi_\al \pi_\be = \de_{\al,\be}\pi_\al$,
$\d\pi_\al=\pi_\al\d$ and comultiplication $\De(\d)=\d\tt 1 + 1\tt\d$,
$\De(\pi_\al) = \sum_{\ga\in\Ga} \pi_{\al\ga^{-1}} \tt \pi_\ga$.  Then
one gets the notion of $\Ga$-twisted conformal algebra \cite{BK} (cf.\
\cite{K2}).

7. Let $D=U(\h)$ be the universal enveloping algebra of the 
Heisenberg Lie algebra $\h$ with generators $a_i$, $b_i$, $c$
and the only non-zero commutation relations
$[a_i, b_i] = c$ ($1\le i\le s$). 
Let $A=DL$ be a free left $D$-module of rank one.
Define $\mu\in P_2(\{A,A\},A)$ by the formula
\begin{equation*}
\mu(L\boxtimes L) = 
\Bigl(\sum_{i=1}^s (a_i\tt b_i - b_i\tt a_i) + c\tt 1 - 1\tt c \Bigr)
\tt_D L.
\end{equation*}
Then $(A,\mu)$ is a Lie algebra in the category $\M^*(D)$ with
annihilation algebra ${K_r}_-$, $r=2s+1$, cf.\ Example \ref{exsr}.

%%%%%%%%%%%%%%%%%%%%%%%%%%%%%%%%%%%%%%%%%%%%%%%%%%%%%%%%%%%%%%%%%%%%%%%%%
\section{Open problems}\label{sop}
%%%%%%%%%%%%%%%%%%%%%%%%%%%%%%%%%%%%%%%%%%%%%%%%%%%%%%%%%%%%%%%%%%%%%%%%%

There are a number of interesting problems which we left beyond the
scope of this paper.

\begin{enumerate}

\item
\begin{sloppypar}
Compute the \coh\ of $\Cur \gtg$ with coefficients in $\Chom (M,N)$,
where $M$ and $N$ are current modules. The same for the Virasoro
conformal algebra, where $M$ and $N$ are modules of densities. Only
$\H^1$ is known (see \cite{CKW}), and the result is highly nontrivial.
\end{sloppypar}

\item
\label{2}
Compute the \coh\ of the general conformal algebra
$\operatorname{gc}_N$ and its infinite-rank subalgebras, see
\cite{K2}, with trivial coefficients. Is it true that $\H^\bullet
(\operatorname{gc}_N, \nc [\d]^N)$ is trivial?

\item
Study the relationship between $\H^\bullet (A,M)$ and $\H^\bullet
(\Lie A, V(M))$. A mapping between the two 
 is given in Section~\ref{relation}. Our
computations show that in the case of a current or the Virasoro
conformal algebra $A$, the image of
$\H^\bullet(A,\nc)$ contains all generators of
$\H^\bullet (\Lie A, \nc)$.

\item
Compute the \coh\ of conformal algebras in several indeterminates.

%\item
%Interpret the diagonal complex of Gelfand--Fuchs 
%(see, \emph{e.g.}, \cite[\S~2.4.1]{F})
%from the point of view
%of conformal \coh.

\item
\begin{sloppypar}
Compute the Hochschild and cyclic conformal \coh\ of $\Cend(M)$. These
problems are apparently related to \ref{2}.
\end{sloppypar}

\end{enumerate}

%%%%%%%%%%%%%%%%%%%%%%%%%%%%%%%%%%%%%%%%%%%%%%%%%%%%%%%%%%%%%%%%%%%%%%%%%
%\bibliographystyle{amsplain}

\end{document}